\documentclass[12pt,reqno]{amsart}

\numberwithin{equation}{section}

\usepackage[final, 
color]{showkeys}
\definecolor{refkey}{gray}{.85}
\definecolor{labelkey}{gray}{.85}

\usepackage{AKstyle}

\usepackage{stmaryrd}
\newcommand{\lp}{\boldsymbol (}
\newcommand{\rp}{\boldsymbol )}

\usepackage{caption}
\usepackage[labelfont=rm]{subcaption}

\usepackage[pagebackref=true, colorlinks]{hyperref}

\hypersetup{pdffitwindow=true,linkcolor=blue,citecolor=blue,urlcolor=blue,menucolor=blue}

\usepackage{comment}

\usepackage[T2A]{fontenc}
\newcommand{\vr}{\varrho}
\newcommand{\Rla}{\mbox{\usefont{T2A}{\rmdefault}{m}{n}\cyrl}}
\renewcommand{\C}{\mathbb C}

\begin{document}

\author{Alex Kontorovich}
\thanks{The author is partially supported by an NSF CAREER grant DMS-1254788 and  DMS-1455705, an NSF FRG grant DMS-1463940, an Alfred P. Sloan Research Fellowship, and a BSF grant.}
\email{alex.kontorovich@rutgers.edu}
\address{Department of Mathematics, Rutgers University, New Brunswick, NJ}

\title[Integral Soddy Sphere Packings]{The Local-Global Principle 
for Integral Soddy Sphere Packings
 }

\begin{abstract}
Fix an integral Soddy sphere packing $\sP$. 
Let $\sB$ be the set of all bends in $\sP$.
A number $n$ is called {\it represented} if $n\in\sB$, that is, if there is a sphere in $\sP$ with bend equal to $n$.
A number $n$ is called {\it admissible} if it is everywhere locally represented, 
meaning
that $n\in\sB(\mod q)$ for all $q$.
It is shown that every sufficiently large admissible number is represented.
\end{abstract}
\date{\today}
\maketitle
\tableofcontents

\newpage

\section{Introduction}


 \begin{figure}
        \begin{subfigure}[t]{0.22\textwidth}
                \centering
\includegraphics[width=\textwidth]{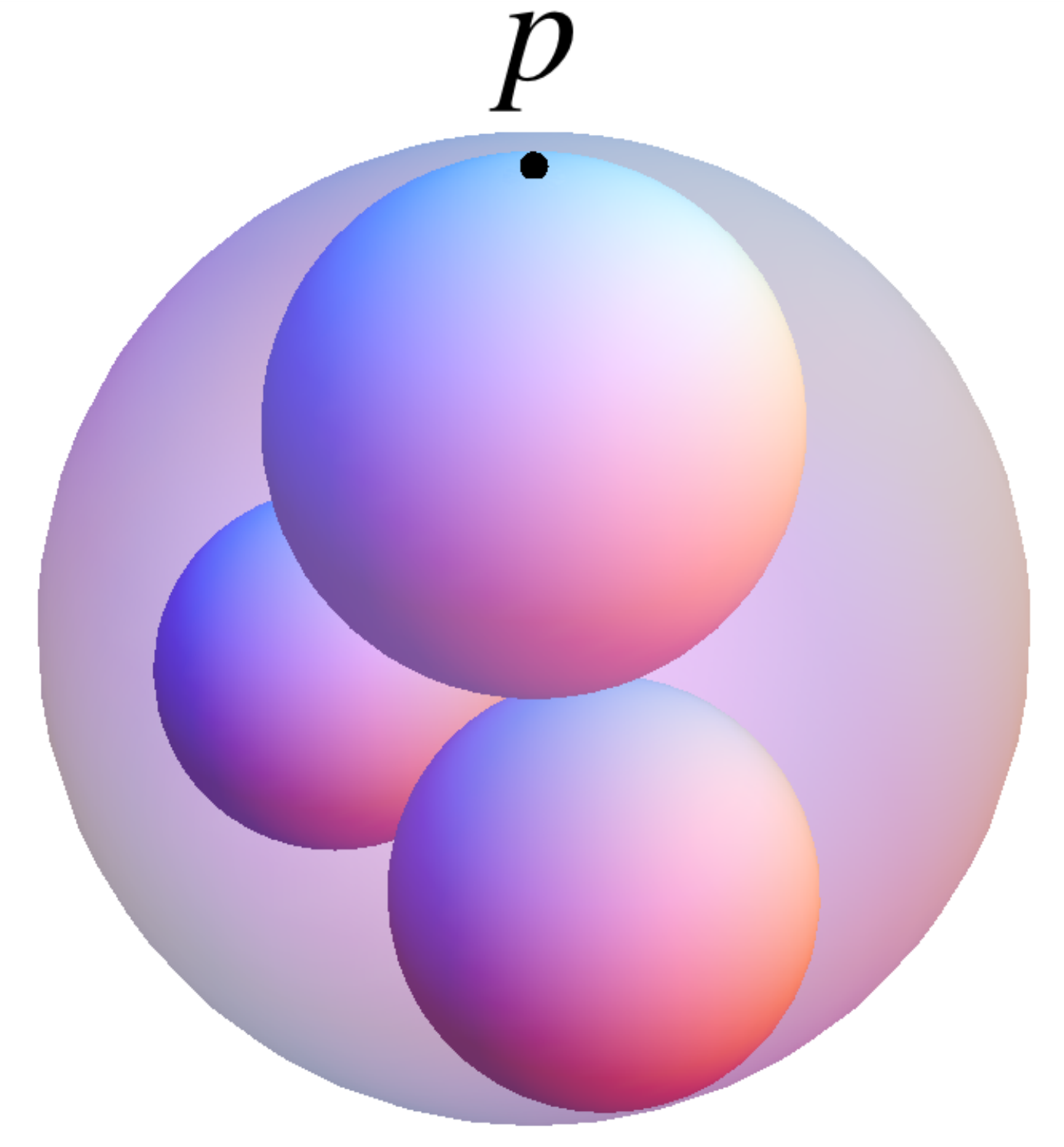}
                \caption{Four tangent spheres}
                \label{fig:Sod4}
        \end{subfigure}%
\quad
        \begin{subfigure}[t]{0.22\textwidth}
                \centering
\includegraphics[width=\textwidth]{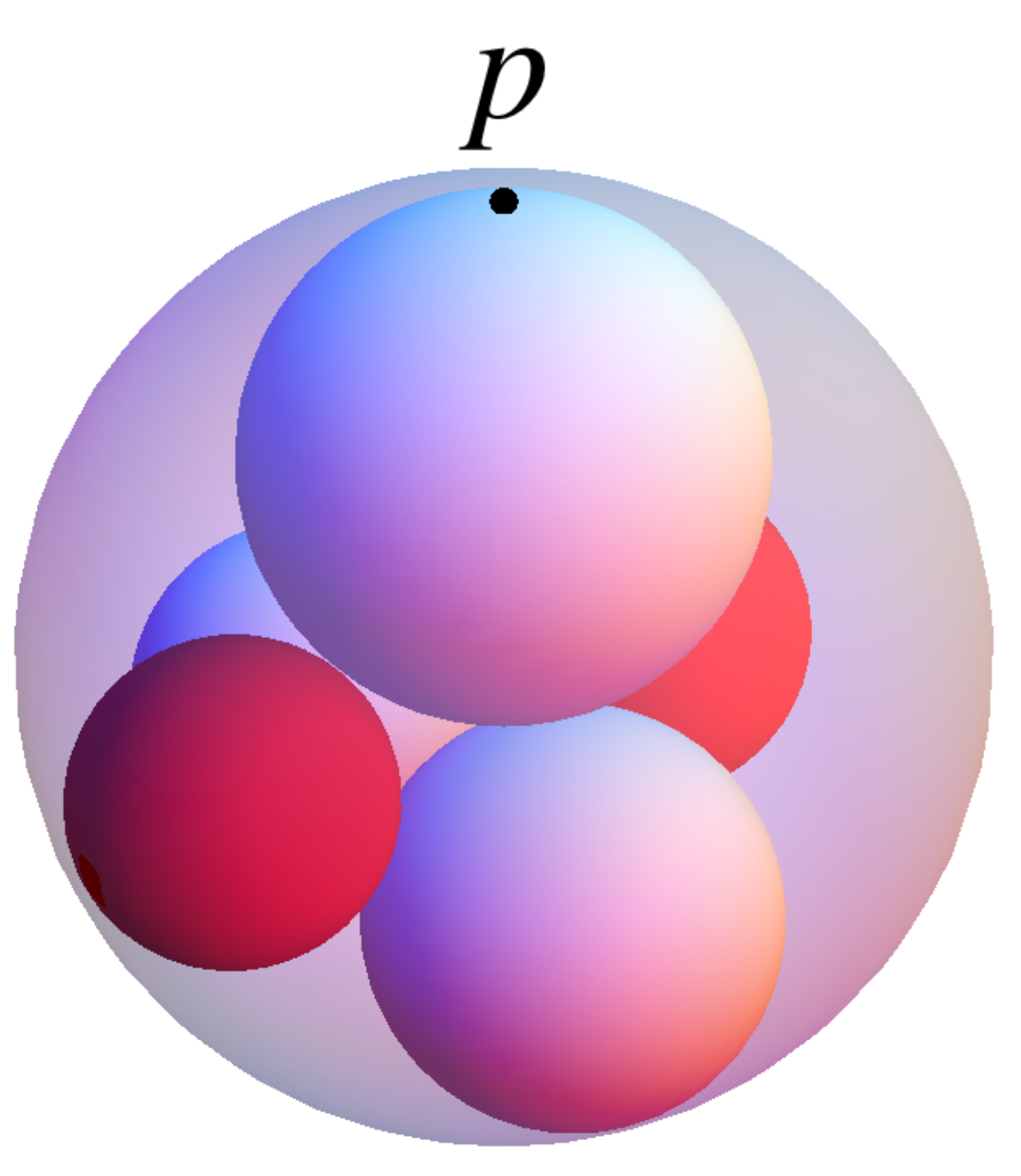}
                \caption{Two more spheres}
                \label{fig:Sod6}
        \end{subfigure}
\quad
        \begin{subfigure}[t]{0.45\textwidth}
                \centering
\includegraphics[width=\textwidth]{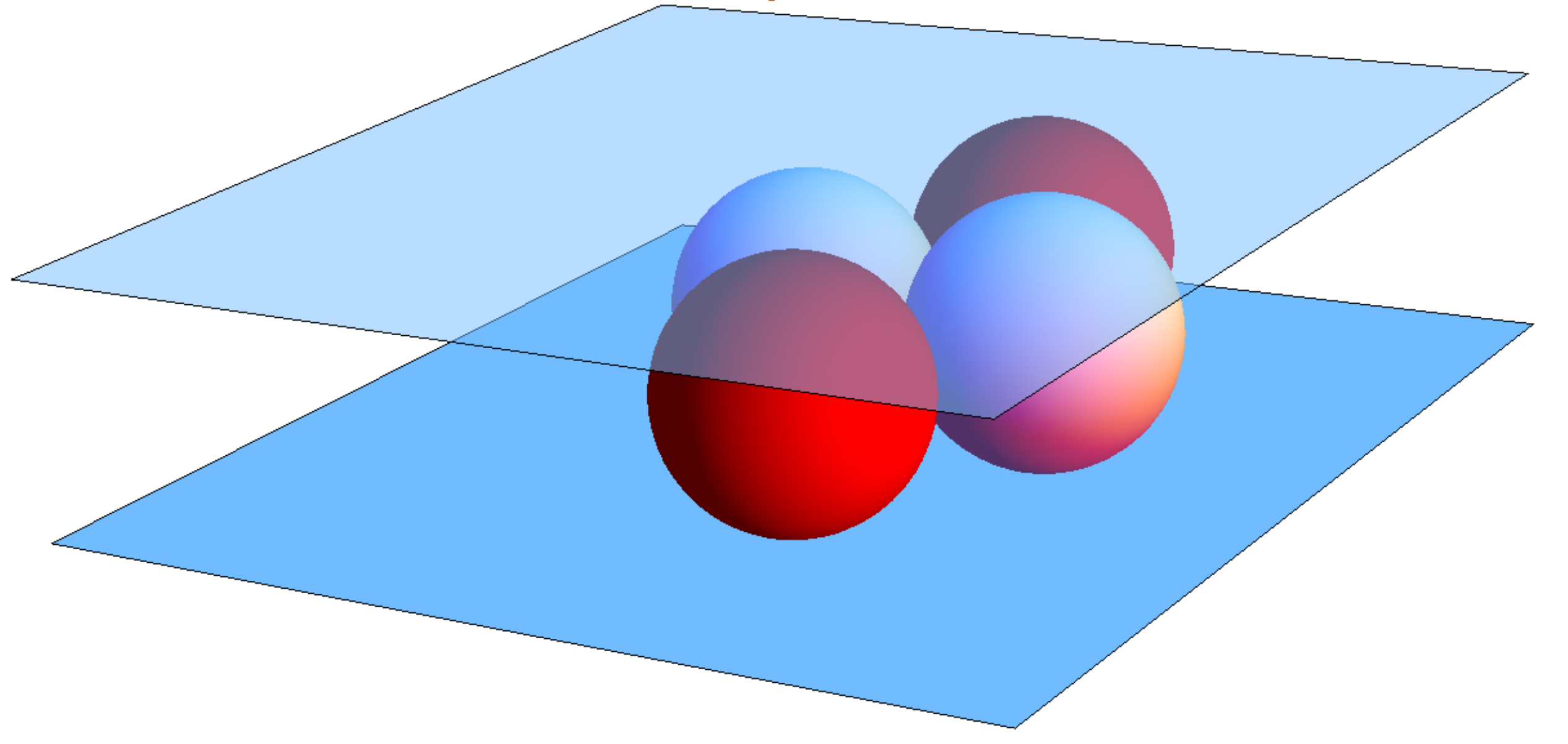}
                \caption{Reflection of
                 (b)
                 through a sphere centered at $p$}
                \label{fig:SodStandard}
        \end{subfigure}%
        \caption{}
\end{figure}

This paper is concerned with
a 
3-dimensional
analogue
 of 
 an
 Apollonian circle packing in the plane, 
 constructed as follows. Given four mutually tangent spheres
  with disjoint points of tangency (\figref{fig:Sod4}),
a generalization to spheres of Apollonius's theorem says that
\be\label{eq:SoddyThm}
\text{
there are exactly two spheres 
}
\ee
tangent to the given ones (\figref{fig:Sod6}).
For a proof of \eqref{eq:SoddyThm}, take a point $p$ of tangency of two given spheres and reflect 
the configuration
through a sphere centered at $p$. 
Thus
$p$ is
 sent to $\infty$,
 and the resulting configuration (\figref{fig:SodStandard}) 
 consists of
 two tangent spheres wedged between two parallel planes; whence the two solutions claimed in \eqref{eq:SoddyThm} are obvious.

\begin{figure}
        \begin{subfigure}[t]{0.41\textwidth}
                \centering
\includegraphics[width=\textwidth]{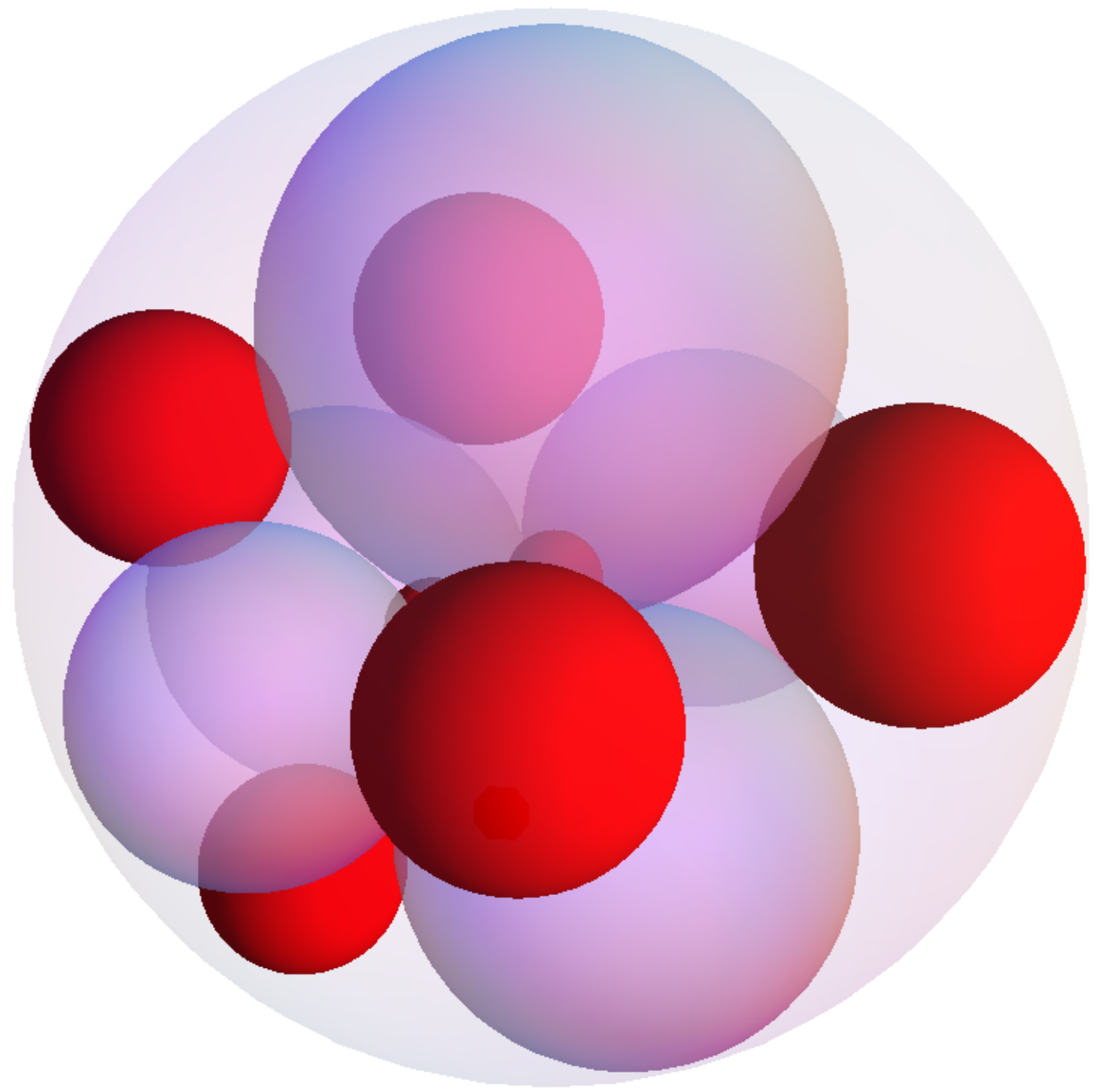}
                \caption{More tangent spheres}
                \label{fig:Sod12}
        \end{subfigure}%
\quad
        \begin{subfigure}[t]{0.52\textwidth}
                \centering
\includegraphics[width=\textwidth]{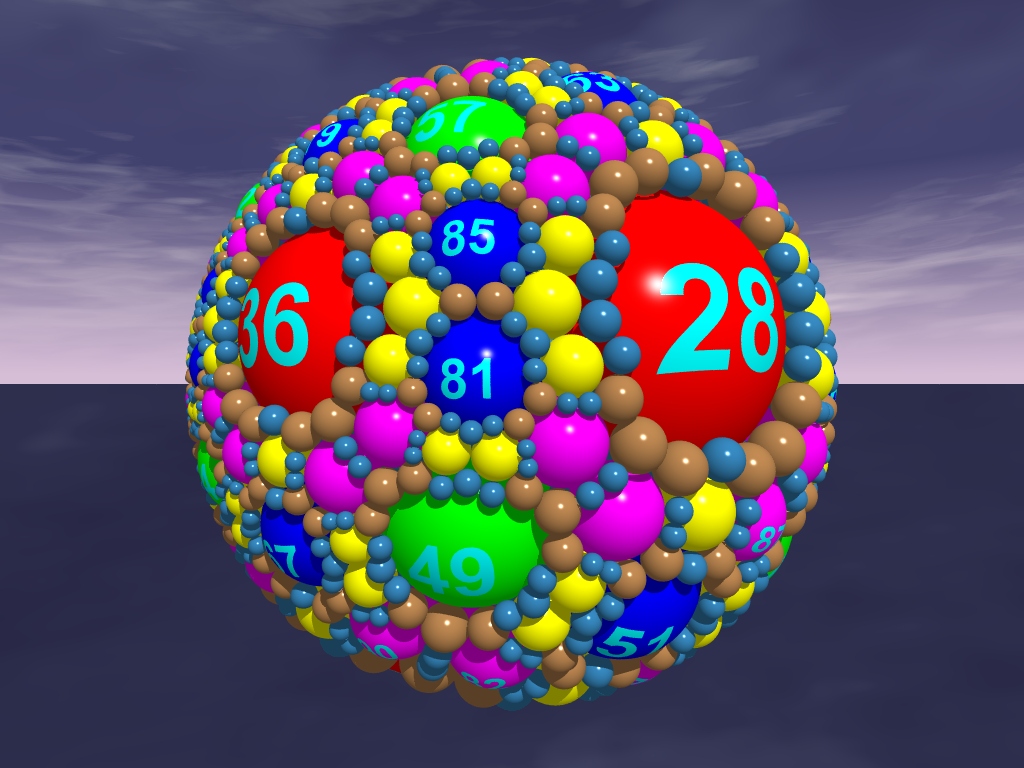}
                \caption{A Soddy sphere packing. Image by\\ Nicolas Hannachi, used with permission}
                \label{fig:Leys}
        \end{subfigure}%
        \caption{}
\end{figure}

Returning to \figref{fig:Sod6}, one 
now
has more configurations of 
tangent spheres, and can iteratively inscribe
further spheres 
in 
this way
 (\figref{fig:Sod12}). Repeating this procedure {\it ad infinitum},
one obtains 
what we will 
call
a {\it Soddy sphere packing} (\figref{fig:Leys}).
\\

\begin{figure}
\includegraphics[width=.4\textwidth]{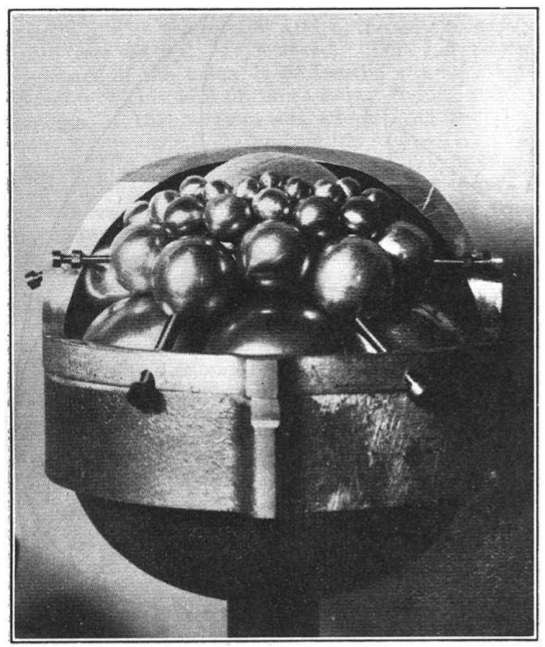}
\caption{A reproduction from \cite{Soddy1937}}
\label{fig:Soddy}
\end{figure}

The name 
 refers 
 to
the
radiochemist
Frederick Soddy (1877-1956),
who in 1936 wrote a {\it Nature} poem
 \cite{Soddy1936} in which he
rediscovered
Descartes's Circle Theorem \cite[pp. 37--50]{Descartes1901} 
and a generalization to spheres, see \thmref{thm:Soddy}.
The latter was known already in 1886 to Lachlan \cite{Lachlan1886}, and  appears in some form as early as 1798 in Japanese Sangaku problems \cite{SangakuW}.
%
%
We 
name the packings after
Soddy because he was the  first, as far as we know, to observe 
that 
there are configurations of circle and sphere packings in which all bends\footnote{The ``bend'' of a circle or sphere is defined to be one over its radius.} are integers \cite{Soddy1937}
; such a packing is called {\it integral}. 
The 
numbers illustrated in \figref{fig:Leys} are some of the bends in that packing. 
In \cite[p. 78]{Soddy1937}, Soddy writes that he ``discovered this [integrality] years ago for the simpler case of cylinders, or circles, in connection with the design of an actual mechanism,'' and provides a picture of a corresponding spherical mechanism, reproduced in \figref{fig:Soddy}.

By rescaling an integral packing, we may assume that the only
 integers
  dividing all of the bends
  are $\pm1$; such a packing is called {\it primitive}.
We restrict 
our attention henceforth to bounded, integral, primitive Soddy sphere packings.
In fact, all of the salient features persist if 
one
considers just
the packing $\sP_{0}$ illustrated in \figref{fig:Leys}.
\\

The goal of this paper is to address the question:
What numbers appear  in \figref{fig:Leys}?
For a 
sphere $S\in\sP$, let $b(S)$ be its bend, and 
let $\sB=\sB(\sP)$ be the set of all bends in $\sP$,
$$
\sB:=\{n\in\Z:\exists S\in\sP,\ b(S)=n\}
.
$$
The bounding sphere is internally tangent to the others, so is given opposite orientation and negative bend.
The first few bends in $\sP_{0}$ are:
\be\label{eq:sKs}
\begin{array}{c}
\sB=\{
-11, 21, 25, 27, 28, 34, 36, 40, 42, 43, 46, 48, 49, 51, 54, 57, 61, 
63,\\    64,67, 69, 70, 72, 73, 75, 78, 79, 81, 82, 84, 85, 87, 90, 
\dots\}.
\end{array}
\ee
A moment's inspection reveals that every bend in $\sP_{0}$ is
\be\label{eq:cong}
\equiv0\text{ or }1\ (\mod3)
,
\ee
that is, there are local obstructions. That such exist was already observed 40 years ago by Boyd \cite[p. 376]{Boyd1973}. 
In analogy with Hilbert's 11th problem on representations of numbers by quadratic forms, we say that $n$ is {\it represented} if $n\in\sB$.
Let $\sA=\sA(\sP)$ be the set of {\it admissible} numbers, that is, numbers $n$ that are everywhere locally represented in the sense that 
\be\label{eq:loc}
n\in\sB(\mod q)\text{ for all $q$}.
\ee 
In our example, $\sA$ is the set of all numbers satisfying \eqref{eq:cong}. The set of admissible numbers 
for any primitive packing $\sP$
satisfies either \eqref{eq:cong} or 
\be\label{eq:cong2}
\equiv0\text{ or }2\ (\mod 3),
\ee 
see \lemref{lem:loc}.

The number of spheres in $\sP$ with bend at most $N$ (counted with multiplicity) is asymptotically equal to a constant times $N^{\gd}$,
where $\gd$ is the Hausdorff dimension of the 
closure of the
packing
(see  \cite{Kim2011}, which generalizes \cite{KontorovichOh2011} to this setting). 
Soddy  
  packings are rigid (one can be mapped to any other by a  conformal transformation), and so $\gd$ is a universal constant;
it is approximately (see \cite{Boyd1973, Borkovec1994}) 
equal to
$$
\gd\approx 2.4739\dots
.
$$
Hence one expects, on grounds of randomness, that the multiplicity of a given admissible bend up to $N$ is roughly $N^{\gd-1}$, which should be quite large. 
In particular, every sufficiently large admissible should be represented.
The main purpose of this paper is to 
confirm this claim.

\begin{thm}[The Local-Global Theorem]\label{thm:main}
The bends of a fixed primitive, integral Soddy sphere packing $\sP$ satisfy a local-to-global principle. 

That is,
there is an effectively computable $N_{0}=N_{0}(\sP)$ so that, if $n>N_{0}$ and $n$ is admissible, $n\in\sA$, then $n$ is represented, $n\in\sB$.
\end{thm}

Empirical evidence suggests
 (and could be verified with enough computation) that 
$N_{0}(\sP_{0})=330$ suffices.\\

\thmref{thm:main} is the analogue to Soddy sphere packings of the local-global conjecture for integral Apollonian circle packings \cite{GrahamLagarias2003, FuchsSanden2011, BourgainFuchs2011, BourgainKontorovich2014a}. Being in higher dimension puts more variables into play, making the problem much easier.
\\

For the proof, we study 
a certain infinite index subgroup $\G$ of the integral orthogonal group 
preserving a particular quadratic form of signature $(4,1)$.
This group, $\G$,   which we call the {\it Soddy group}, is isomorphic to the group of symmetries of $\sP$;
extended to act on hyperbolic $4$-space, the quotient is an infinite volume hyperbolic $4$-fold.
After a calculation,
we find that $\G$ contains an arithmetic (in fact, {\it congruence!}) Kleinian subgroup $\Xi$. 
A consequence is that the set $\sB$ of bends contains the ``primitive'' values of a shifted quaternary quadratic form (and moreover an infinite family of such). 
After some work, we show that
these
satisfy the Hasse principle, from which
the local-global theorem follows. 

This 
proof is a 
generalization to sphere packings of the following related result
 in $2$-dimensions
  due to Sarnak \cite{SarnakToLagarias}: the bends in an integral, primitive Apollonian circle packing contain the primitive values of a shifted {\it binary} quadratic form. 
  It is in this sense that we have more variables: instead of binary forms, sphere packings contain values of quaternary forms.
  Binary 
  forms  
  represent very few numbers, so despite some recent advances \cite{BourgainFuchs2011, BourgainKontorovich2014a}, the analogous problem in circle packings is currently wide open.
So
our main 
innovation
in this paper
is that, for sphere packings, the arithmetic subgroup orbit already fills out all large admissible numbers.

 In dimension $n\ge4$, 
 one can start with a configuration of $n$ tangent hyperspheres, repeating the above-described generating procedure. Unfortunately this does not give rise to a packing, as the hyperspheres eventually overlap \cite{Boyd1973b}.%
 \footnote{Added in print: See Baragar \cite{Baragar2017} for an alternate construction with non-overlapping hyperspheres.}
  Moreover there are no longer any such configurations in which all bends are integral (they can be $S$-integral, with the set $S$ of localized primes depending on the dimension $n$); this follows from Gossett's \cite{Gossett1937} generalization (also in 
 verse) of Soddy's \thmref{thm:Soddy}  to $n$-space.%
 \footnote{Added in print: 
 The tangency graph of a quintuple of mutually tangent spheres generating a Soddy sphere packing is isomorphic to the 1-skeleton of a 4-dimensional simplex. A further generalization in 3-dimensional sphere packings is to consider configurations coming from the 1-skeleton of a 4-orthoplex; the results here have been extended to this setting independently by Dias \cite{Dias2014} and Nakamura \cite{Nakamura2014}. For the corresponding generalization of classical Apollonian packings, see Zhang \cite{Zhang2013a}.}

\subsection*{Notation}
\

The following
  notation for parentheses is used throughout. We sometimes write $x\equiv y(z)$  for $x\equiv y(\mod z)$. 
We will use
bold parentheses
 $\lp x,y\rp$ for the ideal generated by $x$ and $y$, not to be confused with the $\gcd$, denoted simply by $(x,y)$.
 The indicator function $\bo_{\{X\}}$ is $1$ if $X$ holds and $0$ otherwise.

\subsection*{Acknowledgments}\

The author wishes to express his gratitude  to Dimitri Dias, Jeff Lagarias, Yair Minsky,  Kei Nakamura, Alan Reid, and  Peter Sarnak for enlightening conversations, comments and corrections. Thanks also
 to Stony Brook University, where 
the bulk
of 
this text was completed.
\newpage

\section{Preliminaries}


 \begin{figure}
        \begin{subfigure}[t]{0.38\textwidth}
                \centering
\includegraphics[width=\textwidth]{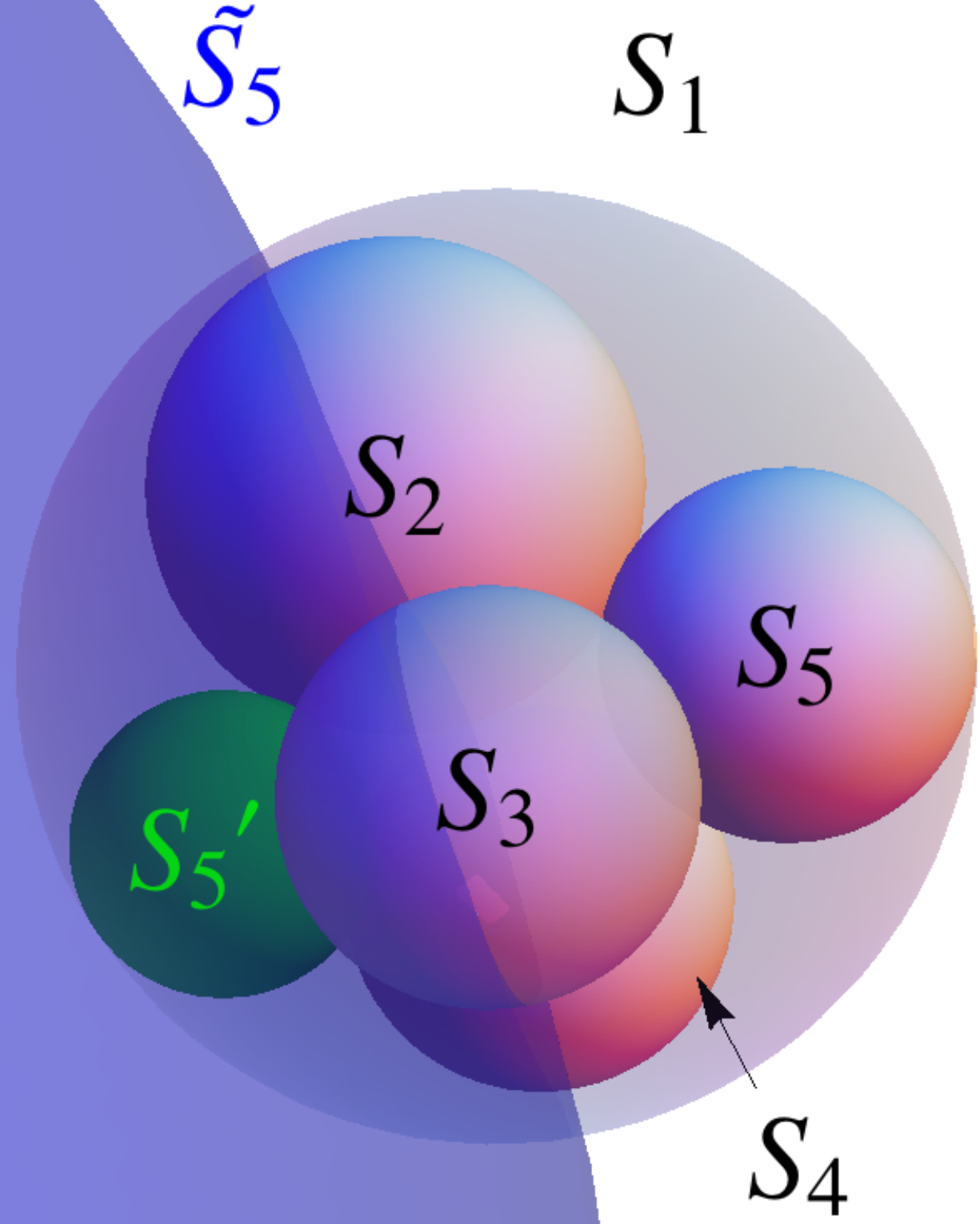}
                \caption{Five tangent spheres 
                and a dual reflection 
                to a sixth
                }
                \label{fig:Dual}
        \end{subfigure}%
\qquad
        \begin{subfigure}[t]{0.53\textwidth}
                \centering
\includegraphics[width=\textwidth]{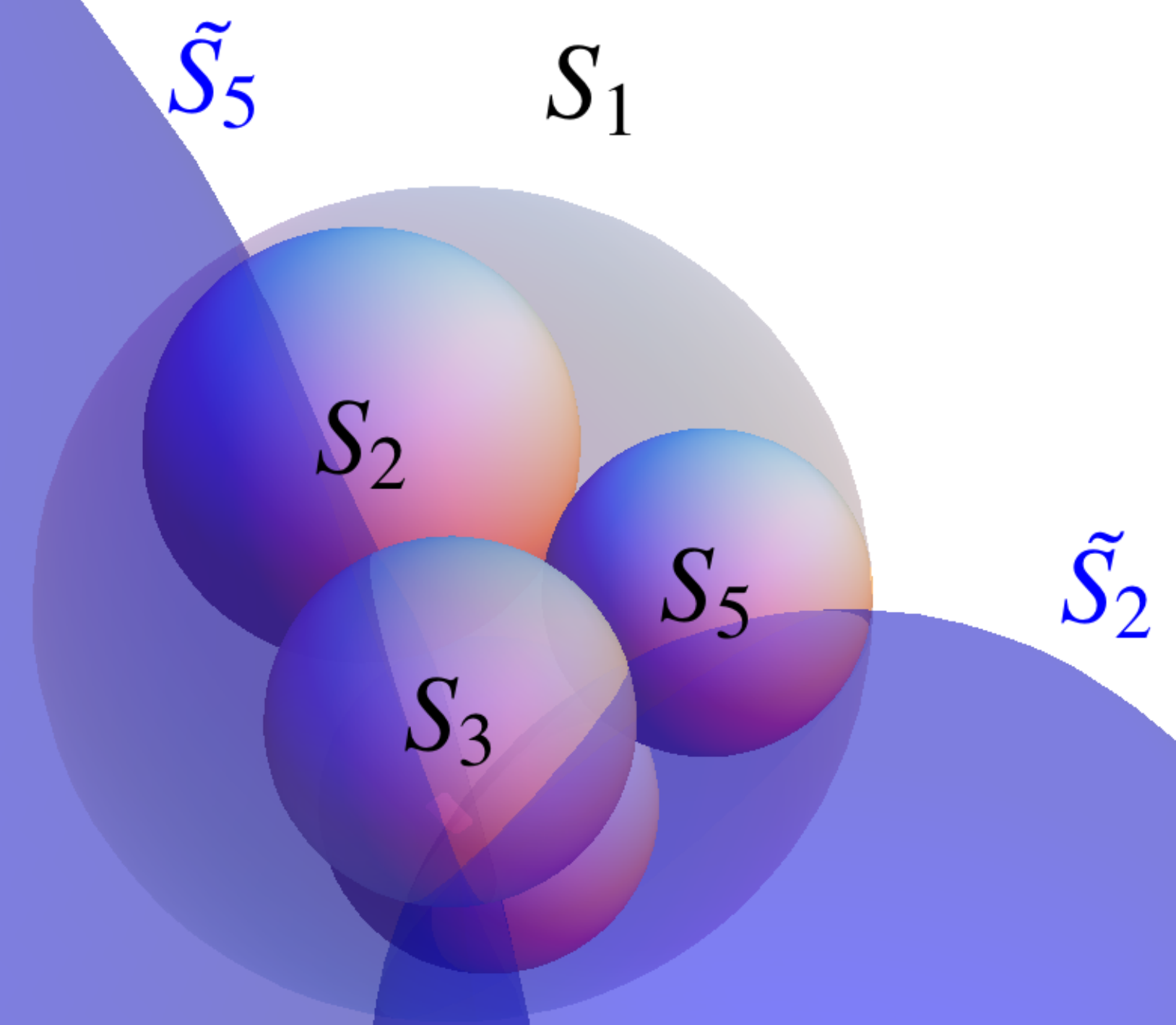}
                \caption{Dual spheres intersect}
                \label{fig:Duals}
        \end{subfigure}%
        \caption{}
\end{figure}

Let 
$\cS=(
S_{1},
S_{2},
S_{3},
S_{4},
S_{5})$ 
 be a configuration of five mutually tangent spheres, and
let 
$$
\bb_{0}=\bb(\cS)
=(
b_{1},
b_{2},
b_{3},
b_{4},
b_{5})
$$ 
 be the corresponding quintuple of bends, with $b_{j}=b(S_{j})$.
Any four tangent spheres, say $S_{1},S_{2},S_{3},S_{4}$ have six cospherical points of tangency, and determine a {\it dual sphere} $\tilde S_{5}$ passing through these points. 
Similarly, for $j=1,\dots,4$, let $\tilde S_{j}$ be the dual sphere orthogonal to all those in $\cS$ except $S_{j}$,
and
call $\tilde\cS=(\tilde S_{1},\dots,\tilde S_{5})$
the {\it dual configuration}.
Reflection 
through $\tilde S_{5}$ fixes $S_{1},S_{2},S_{3},S_{4}$, and sends $S_{5}$ to $S_{5}'$, the other
sphere satisfying \eqref{eq:SoddyThm}, see \figref{fig:Dual}. 
The same holds for the other $\tilde S_{j}$, and iteratively reflecting the original configuration
through the $\tilde S_{j}$ {\it ad infinitum} yields the Soddy packing $\sP=\sP(\cS)$ corresponding to $\cS$.
Observe that unlike the Apollonian case, the dual spheres
 in $\tilde\cS$
 are not tangent, but intersect non-trivially, see \figref{fig:Duals}.

Extend 
the 
reflections through dual spheres
to  hyperbolic $4$-space, 
\be\label{eq:bH4}
\sH^{4}:=\{(x_{1},x_{2},x_{3},y):x_{1},x_{2},x_{3}\in\R,y>0\}
,
\ee
replacing the action of the dual sphere $\tilde S_{j}$ by a reflection
 through a hyper(hemi)sphere
 $\fs_{j}$ 
  whose equator (at $y=0$) is $\tilde S_{j}$ (with $j=1,\dots,5$). We abuse notation, writing $\fs_{j}$ for both the hypersphere and the 
conformal map reflecting through $\fs_{j}$.  
The 
group
\be\label{eq:Gapp}
\cA:=\<
\fs_{1},
\fs_{2},
\fs_{3},
\fs_{4},
\fs_{5}
\>
<
\Isom(\sH^{4})
,
\ee
generated by these reflections 
acts 
discretely on $\sH^{4}$.
 The $\cA$-orbit of any 
 given
 base point in $
 \sH^{4}$ has a limit set in the boundary $\dd\sH^{4}\cong\R^{3}\cup\{\infty\}$, which is  the closure of the original sphere packing. A fundamental domain for this action is the exterior in $\sH^{4}$ of the
 five dual 
hyperspheres $\fs_{j}$.
Hence the quotient hyperbolic $4$-fold $\cA\bk\sH^{4}$ is geometrically finite (with orbifold singularities corresponding to non-trivial intersections of the dual spheres $\tilde S_{j}$), and has infinite hyperbolic volume with respect to 
 the hyperbolic 
 measure
 $$
 y^{-4}dx_{1}dx_{2}dx_{3}dy
$$ 
 in the coordinates \eqref{eq:bH4}.
The group $\cA$ is
the symmetry group of all conformal transformations fixing $\sP$.
\\

For an algebraic realization of $\cA$, we need the following
\begin{thm}[\cite{Lachlan1886, Soddy1936}]\label{thm:Soddy}
Given
a configuration $\cS$ of five tangent spheres, the quintuple $\bb=(b_{1},
b_{2},
b_{3},
b_{4},
b_{5})$ of their bends
lies on the cone 
\be\label{eq:cone}
Q(\bb)=0
,
\ee 
where $Q$ is the quinternary 
quadratic form
\be\label{eq:Qis}
Q(b_{1},\dots,b_{5})
:=
3(
b_{1}^{2}+\cdots+b_{5}^{2}
)-(
b_{1}+\cdots+b_{5}
)^{2}.
\ee
\end{thm}

Recall again that a bounding sphere was negative bend. 
Arguably the nicest formulation of \thmref{thm:Soddy} is the last line of the following excerpt
from Soddy's aforementioned poem \cite{Soddy1936}.
\vskip10pt

{\footnotesize
\qquad To spy out spherical affairs /  
%
An oscular surveyor  /

\qquad Might find the task laborious,  /
%
The sphere is much the gayer,  /

\qquad And now besides the pair of pairs  /
%
A fifth sphere in the kissing shares.  /

\qquad Yet, signs and zero as before,  /
%
For each to kiss the other four  /

\qquad {\it The square of the sum of all five bends  /
%
Is thrice the sum of their squares.}
}
\vskip10pt


If $b_{1},\dots,b_{4}$ are given,
it then
 follows from \eqref{eq:cone} that the variable $b_{5}$ 
 satisfies
 a quadratic equation, and hence there are two solutions.
This is an algebraic proof of \eqref{eq:SoddyThm}. Writing $b_{5}$ and $b_{5}'$ for the two solutions, 
it is elementary from \eqref{eq:cone} that
$$
b_{5}+b_{5}'
=
b_{1}+b_{2}+b_{3}+b_{4}
.
$$
In other words, if the quintuple $(b_{1},b_{2},b_{3},b_{4},b_{5})$ is given, then one obtains the quintuple 
with $b_{5}$ replaced by $b_{5}'$ via a linear action:
$$
\bp
1&&&&\\
&1&&&\\
&&1&&\\
&&&1&\\
1&1&1&1&-1
\ep
\cdot
\bp
b_{1}\\
b_{2}\\
b_{3}\\
b_{4}\\
b_{5}
\ep
=
\bp
b_{1}\\
b_{2}\\
b_{3}\\
b_{4}\\
b_{5}'
\ep
.
$$
This is an algebraic realization of the geometric action of $\tilde S_{5}$ (or $\fs_{5}$) on 
a
quintuple.
Call the above $5\times5$ matrix $M_{5}$. 
One can similarly replace other $b_{j}$ by $b_{j}'$ keeping the four complementary bends fixed, via the matrices
\be\label{eq:Mj}
M_{1}
=
\bp
-1&1&1&1&1\\
&1&&&\\
&&1&&\\
&&&1&\\
&&&&1
\ep
,
M_{2}
=
\bp
1&&&&\\
1&-1&1&1&1\\
&&1&&\\
&&&1&\\
&&&&1
\ep
,
\ee
$$
M_{3}
=
\bp
1&&&&\\
&1&&&\\
1&1&-1&1&1\\
&&&1&\\
&&&&1
\ep
,
M_{4}
=
\bp
1&&&&\\
&1&&&\\
&&1&&\\
1&1&1&-1&1\\
&&&&1
\ep
.
$$
Let
$\G$
be the group  generated by the $M_{j}$:
\be\label{eq:GamAp}
\G:=
\<
M_{1},
M_{2},
M_{3},
M_{4},
M_{5}
\>
.
\ee
By
 construction,
  each generator $M_{j}$ (and hence also $\G$) lies inside 
 the  orthogonal group $O_{Q}$ preserving the 
 form $Q$,
$$
O_{Q}
:=
\left\{
g\in\GL_{5}
:
Q(g\cdot \bb)
=
Q(\bb),\
\forall\bb
\right\}
.
$$
Moreover the Soddy group $\G$ is clearly contained in the group $O_{Q}(\Z)$ of integer matrices.
The fact that $\cA$ has infinite co-volume is equivalent to
 $\G$  having infinite index in $O_{Q}(\Z)$. That is, $\G$ is a ``thin'' group. The generators of $\G$ satisfy the relations: $M_{j}^{2}=I$ and $(M_{j}M_{k})^{3}=I$ \cite[Theorem 5.1]{GrahamLagarias2006}. Geometrically, these relations correspond, respectively, to reflections being involutions, and to the non-trivial intersections of the dual spheres (recall \figref{fig:Duals}).

The orbit 
\be\label{eq:cOAp}
\sO:=\G\cdot\bb
\ee
of the quintuple $\bb
=\bb(\cS)$ under the Soddy group $\G$ 
consists of all quintuples corresponding to bends of five mutually tangent spheres in the packing $\sP$.
Hence the set $\sB$ of all bends in $\sP$ is simply the union of sets of the form
\be\label{eq:KwGv}
\sB=\bigcup_{\bw
\in\{\bbe_{1},\dots,\bbe_{5}\}}\<\bw
,\G\cdot\bb
\>
,
\ee
as $\bw
$ ranges through the standard basis vectors 
$$
\bbe_{1}
=
(1,0,0,0,0), 
 \cdots ,
  \bbe_{5}
  =
  (0,0,0,0,1).
  $$
  The inner product $\<\cdot,\cdot\>$ in \eqref{eq:KwGv} is the standard one on $\R^{5}$.
  
This 
explains 
the integrality of all bends in \figref{fig:Leys}: the group $\G$ has only integer matrices, so if the initial quintuple $\bb_{0}$ (or for that matter, any 
bends of 
five mutually tangent spheres in $\sP$) is integral, then the bends in $\sP$ are all integers (as first observed by Soddy \cite{Soddy1937}). \\

From \eqref{eq:KwGv} it is elementary to
see the local obstruction claimed in \eqref{eq:cong}. For the packing $\sP_{0}$ of \figref{fig:Leys}, one can choose to
 generate from
  the
``root'' quintuple (meaning it consists of the bends of the five largest tangent spheres, see \cite[\S3]{GrahamLagarias2003}) 
\be\label{eq:bv0Is}
\bb_{0}:=(-11, 21, 25, 27, 28).
\ee
The orbit under $\G$, reduced mod $3$, is then elementarily computed. In general we have the following
\begin{lem}\label{lem:loc}
For $\sB$ the set of bends of an integral, primitive Soddy packing $\sP$, there is always a local obstruction $\mod 3$, either of the form \eqref{eq:cong} or \eqref{eq:cong2}. In particular, there is an $\vep=\vep(\sP)\in\{\pm1\}$ so that, for any quintuple $\bb$ in the cone \eqref{eq:cone} over $\Z$, two entries are $\equiv0(\mod 3)$ and three entries are $\equiv\vep(\mod 3)$.
\end{lem}
Note that we are not (yet) claiming that these are the {\it only} local obstructions; this will follow from our proof of the local-to-global theorem.
\pf
One 
may
first 
attempt to understand
the cone
 \eqref{eq:cone}  over $\Z/3\Z$, but the form $Q$ in \eqref{eq:Qis} reduced  mod $3$ is highly degenerate. So instead consider the cone 
   over $\Z/9\Z$. 
Disregarding   
 the origin (since the packing is assumed to be primitive),
 there are $140$ vectors $\mod 9$, not counting permutations. Reducing these $\mod 3$ leaves only
the two 
 vectors $(0,0,\vep,\vep,\vep)$, $\vep\in\{\pm1\}$, and their permutations. 
The action of $\G(\mod 3)$ on these 
is 
trivial: each vector is fixed. 
This is all verified by direct computation.
\epf

It is convenient to also record here the following

\begin{lem}\label{lem:bend24mod6}
The set  $\sB$ of bends of an integral, primitive Soddy packing $\sP$ always contains 
an element $b\equiv\vep(\mod 6)$, and
an element 
$b\equiv3+\vep(\mod 6)$.
\end{lem}
\pf
The cone \eqref{eq:cone} mod $36$ has $30,576$ vectors, not counting permutations. Reducing these mod $6$ leaves $15$ vectors, of which $5$ are imprimitive, the remaining ones being:
$$
\text{if $\vep(\sP)=+1$:}\qquad
\left(
\begin{array}{c}
 0 \\
 0 \\
 1 \\
 1 \\
 1 \\
\end{array}
\right)
,
\left(
\begin{array}{c}
 0 \\
 1 \\
 1 \\
 1 \\
 3 \\
\end{array}
\right)
,\left(
\begin{array}{c}
 0 \\
 1 \\
 1 \\
 3 \\
 4 \\
\end{array}
\right)
,
\left(
\begin{array}{c}
 1 \\
 1 \\
 3 \\
 3 \\
 4 \\
\end{array}
\right)
,
\left(
\begin{array}{c}
 1 \\
 3 \\
 3 \\
 4 \\
 4 \\
\end{array}
\right)
$$
$$
\text{if $\vep(\sP)=-1$:}\qquad
\left(
\begin{array}{c}
 0 \\
 0 \\
 5 \\
 5 \\
 5 \\
\end{array}
\right)
,
\left(
\begin{array}{c}
 0 \\
 3 \\
 5 \\
 5 \\
 5 \\
\end{array}
\right)
,
\left(
\begin{array}{c}
 0 \\
 2 \\
 3 \\
 5 \\
 5 \\
\end{array}
\right)
,
\left(
\begin{array}{c}
 2 \\
 2 \\
 3 \\
 3 \\
 5 \\
\end{array}
\right)
,
\left(
\begin{array}{c}
 2 \\
 3 \\
 3 \\
 5 \\
 5 \\
\end{array}
\right)
.
$$
They plainly each have at least one element $\equiv1$ or $5(\mod 6)$, giving the first claim.

Next observe that the orbit under $\G(\mod 6)$ plus permutations acts transitively on each row (of course, $\G$ cannot change $\vep(\sP)$). This gives the second claim, that one can always make either $2$ or $4$ appear as one of the entries $\mod 6$.
\epf

\newpage

\section{Bends as Primitive Values of Quaternary Forms}

In this section, we show that a  subset of the  bends $\sB$ in a Soddy packing can be obtained as ``primitive'' 
(which has a  non-standard meaning here; see below)
values of certain shifted quaternary quadratic forms.
%
%
Our first goal is to
prove that the Soddy group $\G$, while being infinite index in $O_{Q}\cong O(4,1)$, contains a
congruence Kleinian subgroup.
The method is a generalization of Sarnak's observation in \cite{SarnakToLagarias}.

\begin{figure}
        \begin{subfigure}[t]{
        \textwidth}
                \centering
\includegraphics[width=\textwidth]{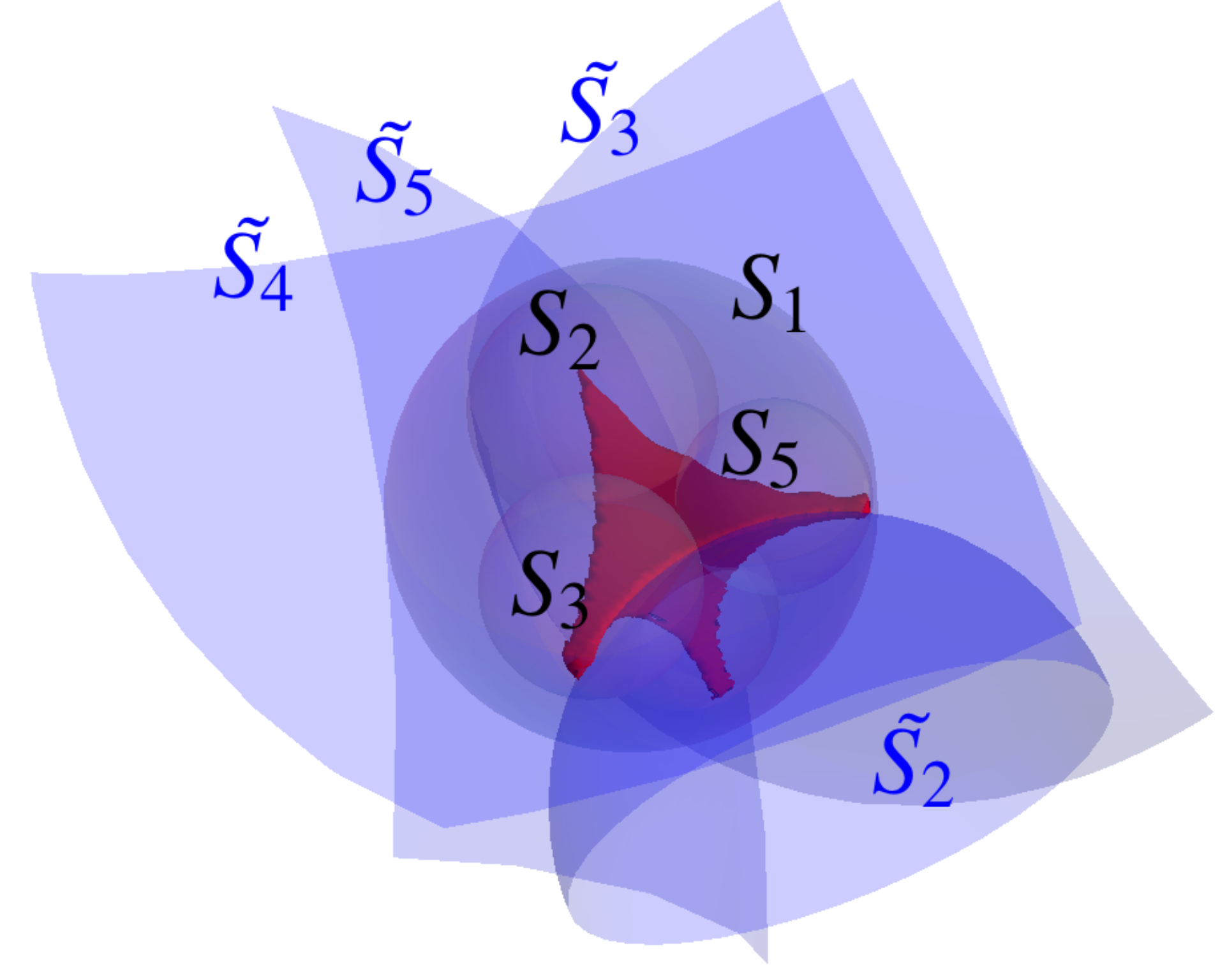}
\caption{A fundamental domain for the action of $\cA_{1}$}
                \label{fig:A1}
        \end{subfigure}%
\qquad
        \begin{subfigure}[t]{
        \textwidth}
                \centering
\includegraphics[width=\textwidth]{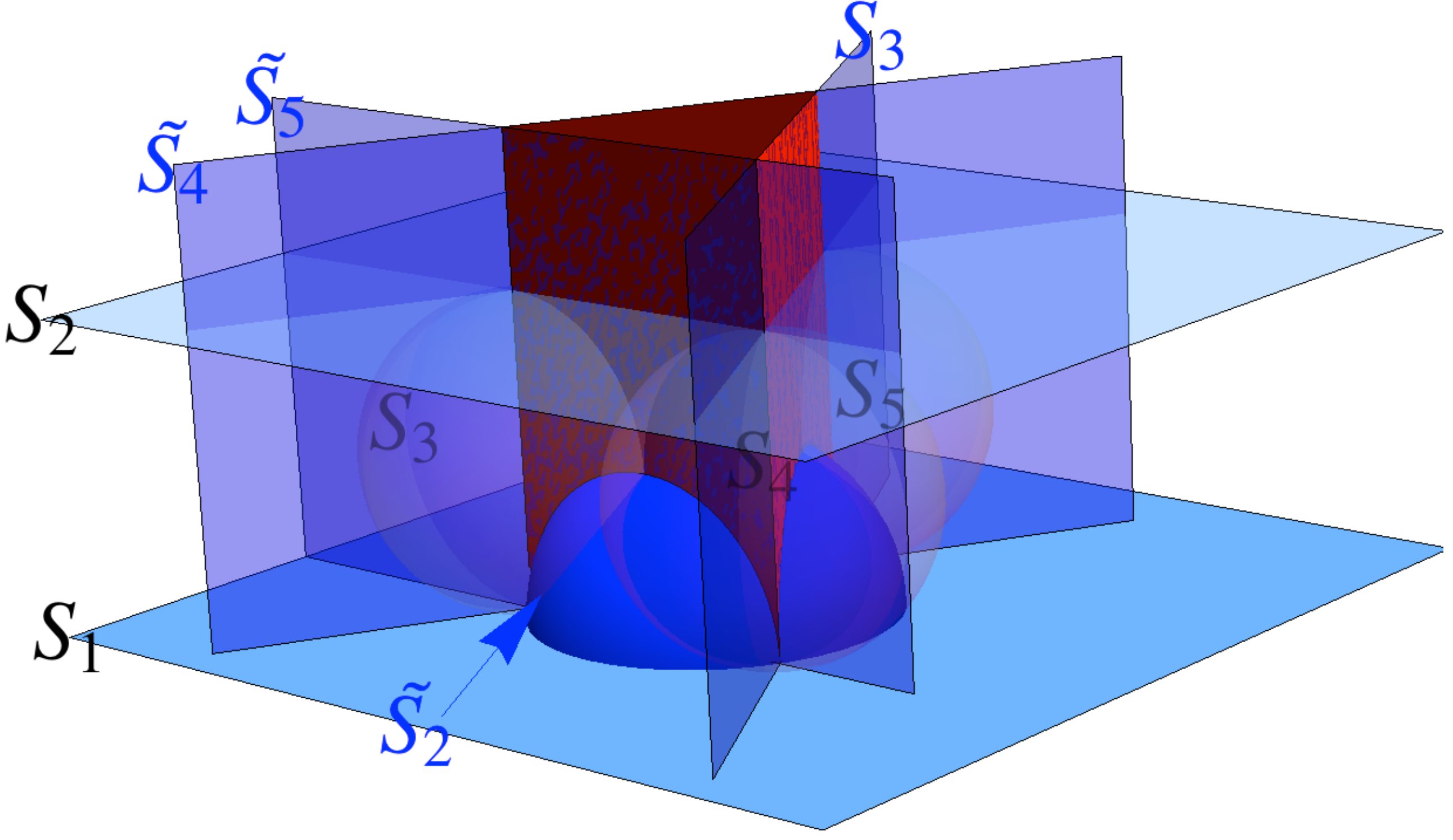}
\caption{The same domain on sending two spheres to planes}
\label{fig:A1std}
        \end{subfigure}%
        \caption{}
\end{figure}

Recall
the configuration $\cS=(S_{1},\dots,S_{5})$ of five mutually tangent spheres 
and
 the group $\cA$ in \eqref{eq:Gapp} 
of reflections through spheres in the configuration $\tilde\cS$ dual to $\cS$.
Let 
$$
\cA_{1}=\<\fs_{2},...,\fs_{5}\>
$$ 
be the subgroup of $\cA$ which fixes the sphere $S_{1}$ in $\cS$.
It acts discontinuously on the interior of $S_{1}$, 
which we now consider as 
the ball model for hyperbolic $3$-space $\sH^{3}$.
A fundamental domain for 
the quotient $\cA_{1}\bk\sH^{3}$
 is
  the
 curvilinear
 regular
 ideal  tetrahedron interior to $S_{1}$ and exterior to the dual spheres $\tilde S_{2},\dots,\tilde S_{5}$, see \figref{fig:A1}.
This is easier seen by first 
applying the same
 transformation as in \figref{fig:SodStandard}, see \figref{fig:A1std}.
   In particular, the quotient has  finite volume, 
   and
at    any vertex,
   the three edges meet at
 dihedral angles 
 all equal to $\pi/3$.   Then the volume
   can be computed via the dilogarithm, or equivalently, Lobachevsky's function
$$
 \Rla(\gt):=-\int_{0}^{\gt}\log|2\sin u|\, du
 ,
$$
  see, e.g., \cite[Lemma 2]{Milnor1982}.
Namely, the volume of this domain is 
$
 3\,\Rla(\pi/3).
$ 
Then its index-2 orientation preserving subgroup, a gluing of two such tetrahedra, has co-volume
 \be\label{eq:volA1}
 \vol((\cA_{1}\cap\Isom^{+})\bk\sH^{3})=6\,\Rla(\pi/3).
 \ee 

\begin{rmk}
\label{rmk:McM}
Curt McMullen asked (private communication) whether  this 
quotient
is then the figure eight knot complement; recall that Thurston 
 showed 
the latter
can be triangulated by two maximal tetrahedra. 
It turns out that,
like the knot complement,
 our quotient is indeed arithmetic;
but 
the two are not isomorphic, see \rmkref{rmk:Grun} below.
\end{rmk}

To realize this geometric action algebraically, let 
\be\label{eq:G1is}
\G_{1}:=\<M_{2},\cdots,M_{5}\>
\ee
be the corresponding subgroup of $\G$, where the $M_{j}$ are given in \eqref{eq:Mj}.
We immediately pass again to the 
index-2
orientation preserving subgroup, setting  
\be\label{eq:XiIs}
\Xi:=\G_{1}\cap\SL_{5}.
\ee
Then $\Xi$ is generated by
\be\label{eq:Xidef}
\Xi=
\<
\xi_{1},\xi_{2},\xi_{3}
\>
,
\ee
where
$$
\xi_{1}:=M_{2}M_{3}=\left(
\begin{array}{ccccc}
 1 & 0 & 0 & 0 & 0 \\
 2 & 0 & -1 & 2 & 2 \\
 1 & 1 & -1 & 1 & 1 \\
 0 & 0 & 0 & 1 & 0 \\
 0 & 0 & 0 & 0 & 1
\end{array}
\right),\ 
\xi_{2}:=
M_{2}M_{4}
=
\left(
\begin{array}{ccccc}
 1 & 0 & 0 & 0 & 0 \\
 2 & 0 & 2 & -1 & 2 \\
 0 & 0 & 1 & 0 & 0 \\
 1 & 1 & 1 & -1 & 1 \\
 0 & 0 & 0 & 0 & 1
\end{array}
\right)
,
$$
and
$$
\xi_{3}:=
M_{2}M_{5}=
\left(
\begin{array}{ccccc}
 1 & 0 & 0 & 0 & 0 \\
 2 & 0 & 2 & 2 & -1 \\
 0 & 0 & 1 & 0 & 0 \\
 0 & 0 & 0 & 1 & 0 \\
 1 & 1 & 1 & 1 & -1
\end{array}
\right)
.
$$

It will turn out that $\Xi$ is in fact a {\it congruence} group, as a form of $\SL_{2}(\C)$. To see this, we make a number of transformations. 

\begin{lem}
Let 
$$
J
=
\left(
\begin{array}{ccccc}
 1 & 0 & 0 & 0 & 0 \\
 1 & 0 & 0 & 1 & 0 \\
 1/3 & 1/3 & -2/3 & 1/3 & 1/3 \\
 1/3 & -2/3 & 1/3 & 1/3 & 1/3 \\
 1 & 0 & 0 & 0 & 1
\end{array}
\right)
.
$$
Then for $j=1,2,3$, the conjugates 
\be\label{eq:Jcong}
\tilde\xi_{j}:=J\cdot\xi_{j}\cdot J^{-1}
\ee
 are given by
\be\label{eq:tilXis}
\tilde\xi_{1}=
\left(
\begin{array}{ccccc}
 1 & 0 & 0 & 0 & 0 \\
 0 & 1 & 0 & 0 & 0 \\
 0 & 0 & 0 & 1 & 0 \\
 0 & 0 & -1 & -1 & 0 \\
 0 & 0 & 0 & 0 & 1
\end{array}
\right)
,\
\tilde\xi_{2}=
\left(
\begin{array}{ccccc}
 1 & 0 & 0 & 0 & 0 \\
 0 & 1 & -3 & -3 & 3 \\
 0 & 0 & -1 & -1 & 2 \\
 0 & 0 & 1 & 0 & -1 \\
 0 & 0 & 0 & 0 & 1
\end{array}
\right)
,
\ee
and
$$
\tilde\xi_{3}=
\left(
\begin{array}{ccccc}
 1 & 0 & 0 & 0 & 0 \\
 0 & 1 & 0 & 0 & 0 \\
 0 & 2 & -1 & -1 & 0 \\
 0 & -1 & 1 & 0 & 0 \\
 0 & 3 & -3 & -3 & 1
\end{array}
\right)
.
$$
\end{lem}
\pf
Of course this can be verified by direct computation. But we elucidate the role of $J$ as follows.

Let $\bb=\bb(\cS)=(b_{1},\dots,b_{5})$ be the quintuple of bends corresponding to $\cS$.
Write the 
form $Q$ in \eqref{eq:Qis} as
\beann
Q(b_{1},b_{2},\dots,b_{5})
&=&
3(b_{1}^{2}+b^{2}_{2}+\cdots+b^{2}_{5})-(b_{1}+b_{2}+\cdots+b_{5})^{2}
\\
&=&
2(\widetilde Q(\by)+3b_{1}^{2})
,
\eeann
where
\be\label{eq:by}
\by=(y_{2},\dots,y_{5})
:=(b_{2},\dots,b_{5})+(b_{1},b_{1},b_{1},b_{1}),
\ee
and
$$
\widetilde Q(\by):=
y_{2}^{2}+\cdots+y_{5}^{2}-y_{2}y_{3}-y_{2}y_{4}-\cdots-y_{4}y_{5}.
$$
The affine action of $\Xi$ on $(b_{2},\dots,b_{5})$ is conjugated by the above to a linear action $\Xi'<\SO_{\widetilde Q}$.
Since $\bb$ was assumed to be 
primitive, $\by
$ is a primitive point on the quadric
\be\label{eq:gyIs}
\widetilde Q(\by)=-3b_{1}^{2}.
\ee
%
For
later
convenience, 
 we make another
 change of variables. 
 It turns out that, despite  beginning with a problem in the (rational) integers, we will need to work in the number field
$$
K:=\Q(\sqrt{-3})
$$ 
with its ring of (Eisenstein) integers
$$
\cO:=\Z[\gw].
$$
Here 
$$
\gw:=e^{\pi i/3}
$$
is a primitive {\it sixth} root of unity (it turns out to be more convenient to use the sixth root than the cube root). 
We will
 conjugate $\widetilde Q$ to 
 the form
\be\label{eq:Fis}
F(\bba):=B^{2}+BC+C^{2}-AD
,
\ee
where $\bba=(A,B,C,D)$.
The determinant of the Hermitian matrix
\be\label{eq:Xis}
X:=
\mattwo{A}{B+\gw C}{B+\bar\gw C}{D}
\ee
is easily seen to be $-F(\bba)$.
Let
$$
y_{2}=A - B- 2 C+ D
,\
y_{3}=A - 2 B- C+ D
,\
y_{4}=A
,\
y_{5}=D
,
$$
or equivalently, make the change of variables
$$
A=y_{4}
,\
B=
\frac{y_{2}-2y_{3}+y_{4}+y_{5}}3
,\
{C}=
\frac{-2y_{2}+y_{3}+ y_{4}+ y_{5}}3
,\
{D}=
{y_{5}}
.
$$
We claim  that $B$ and $C$ are integers; indeed, returning to the $b$ variables in \eqref{eq:by}, we have
\bea\nonumber
A&=&  b _1+  b _4 ,\\
\nonumber
B&=&\frac13( b _1+ b _2-2  b _3+ b _4+ b _5)
,\\
\nonumber
C&=&\frac13( b _1-2  b _2+ b _3+ b _4+ b _5)
,\\
\label{eq:ABCD}
D&=& b _1+ b _5
.
\eea
But  
reducing \eqref{eq:cone}, \eqref{eq:Qis} mod $3$ shows that $b_{1}+\dots+b_{5}\equiv0(\mod 3)$, and hence $B$ and $C$ are integers.

In these coordinates, \eqref{eq:gyIs} becomes
\be\label{eq:Ftogk1}
F(\bba)
=-b_{1}^{2}
.
\ee
The action $\Xi'<\SO_{\widetilde Q}$ on $\by$ is then conjugated to an action $\widetilde\Xi<\SO_F$ on $\bba$.

The matrix $J$ is then simply the change of variables matrix from $\bb$ to $(b_{1},\bba)$. 
\epf
The convenience of 
this conjugation is made apparent in the following 
\begin{lem}
The quadratic form $F$ in \eqref{eq:Fis} has signature $(3,1)$. The connected component of the identity of the special orthogonal group $\SO_{F}(\R)$ 
has spin double cover isomorphic to $\PSL_{2}(\C)$.
There is a homomorphism $
\rho:\PSL_{2}(\C)\to\SO_{F}(\R)$ given explicitly (for our purposes embedded in $\GL_{5}$) by
mapping
\be\label{eq:g}
g=\mattwo
\ga\gb\g\gd
\in\PSL_{2}(\C)
\ee
to $\frac1{|\det(g)|^{2}}\times$
\be\label{eq:gTil}
\hskip-20pt
\left(
\begin{array}{ccccc}
1&&&&\\
&
 |\alpha |^2 & 2 \Re\left(\beta  \bar{\alpha }\right) & 2 \Re\left(\alpha  \omega  \bar{\beta
   }\right) & |\beta |^2 \\
   &
 \frac{2 }{\sqrt{3}}\Im\left(\gamma  \omega  \bar{\alpha }\right) & \frac{2}{\sqrt{3}} \Im\left(\omega 
   \left(\delta  \bar{\alpha }+\gamma  \bar{\beta }\right)\right) & \frac{2}{\sqrt{3}}    \Im\left(\gamma  \bar{\beta } \omega ^2+\delta  \bar{\alpha }\right) & \frac{2}{\sqrt{3}}\Im\left(\delta  \omega  \bar{\beta }\right) \\
   &
 \frac{2}{\sqrt{3}} \Im\left(\alpha  \bar{\gamma }\right) & \frac{2}{\sqrt{3}} \Im\left(\beta  \bar{\gamma
   }+\alpha  \bar{\delta }\right) & \frac{2 }{\sqrt{3}}\Im\left(\omega  \left(\alpha  \bar{\delta
   }-\gamma  \bar{\beta }\right)\right) & \frac{2 }{\sqrt{3}}\Im\left(\beta  \bar{\delta
   }\right) \\
   &
 |\gamma |^2 & 2 \Re\left(\gamma  \bar{\delta }\right) & 2 \Re\left(\gamma  \omega  \bar{\delta
   }\right) & |\delta |^2
\end{array}
\right)
.
\ee
The  preimages under $
\rho$ of the matrices $\tilde\xi_{1},\tilde\xi_{2},\tilde\xi_{3}$ in \eqref{eq:tilXis}  are $\pm \ft_{1},\pm\ft_{2},\pm\ft_{3}$, respectively, where:
\be\label{eq:SpinMats}
\ft_{1}=
\left(
\begin{array}{cc}
{\omega^{-1} } & 0 \\
 0 & \omega 
\end{array}
\right)
,\quad
\ft_{2}=
\left(
\begin{array}{cc}
{\omega ^{-2}} & 
\omega\vr  \\
 0 & \omega ^2
\end{array}
\right)
,\quad
\ft_{3}=
\left(
\begin{array}{cc}
 \omega  & 0 \\
\omega\vr  & {\omega^{-1} }
\end{array}
\right)
.
\ee
Here 
$$
\vr:=1+\gw
$$
is the prime in $\cO$ above the ramified rational prime $3$, which factors as
$
3=\bar\gw \vr^{2}.
$
\end{lem}
\pf
The signature of $F$ is computed directly, and its spin group being $\PSL_{2}(\C)$ is a general fact in the theory of quadratic forms, see e.g. \cite[Ch. 10]{Cassels1978}. We construct $
\rho$ explicitly as follows.
Return to the
 Hermitian matrix 
 $X$ in \eqref{eq:Xis} with
 determinant $-F(\bba)$.
Then for $g\in\PSL_{2}(\C)$,
$$
X':= g\cdot X\cdot \bar g^{t} = 
\mattwo{A'}{B'+\gw C'}{B'+\bar\gw C' }{D'}
$$ 
is also Hermitian with the same determinant.
This gives a linear action sending $(A,B,C,D)$ to $(A',B',C',D')$,
which can be computed explicitly in the coordinates \eqref{eq:g}. The result (embedded in $\GL_{5}$)
is \eqref{eq:gTil}. The preimages \eqref{eq:SpinMats} are then computed directly.
\epf

Let 
\be\label{eq:gLis}
\gL=\<\pm\ft_{1},\pm\ft_{2},\pm\ft_{3}\>/\<\pm I\> \ <\  \PSL_{2}(\C)
\ee
be the group generated by \eqref{eq:SpinMats}.

Then $\gL$
is clearly a subgroup 
of the Bianchi group $\PSL_{2}(\cO)$.
The full group  
$\PSL_{2}(\cO)$ is well-known to have co-volume 
$$
\vol(\PSL_{2}(\cO)\bk\sH^{3})=\foh\Rla(\pi/3),
$$ 
see e.g. \cite[p. 21]{Milnor1982}. 
Combined with
\eqref{eq:volA1}, this
gives us the index
$$
[\PSL_{2}(\cO):\gL]=12,
$$
%
since $\gL\cong\Xi\cong\cA_{1}\cap\Isom^{+}$.
\begin{rmk}\label{rmk:Grun}
This fact was already known to Grunewald-Schwermer, who list a conjugate of the generators \eqref{eq:SpinMats} in their table \cite[p. 76]{GrunewaldSchwermer1993}, calling the group ``$\G_{-3}(12,7)$''.
In the same table [p. 75], the figure eight knot complement is listed as ``$\G_{-3}(12,1)$''; so these are not 
isomorphic, cf. \rmkref{rmk:McM}.
\end{rmk}
The next lemma, crucial for our purposes, states that
our group is not just arithmetic, but {\it congruence}.

\begin{lem}
The group $\gL$ is 
equal to
the following 
congruence 
subgroup 
of $\PSL_{2}(\cO)$,
\be\label{eq:Gam0Is}
\left\{
\mattwo\ga\gb\g\gd\in\PSL_{2}(\cO):\gb,\g\equiv0(\mod \vr)
\right\}
.
\ee

\end{lem}

\pf
The inclusion 
\be\label{eq:gLinG}
\gL\quad<\quad
\eqref{eq:Gam0Is}
\ee 
is clear from the 
generators \eqref{eq:SpinMats}. For the opposite inclusion,
is it an elementary computation that
\eqref{eq:Gam0Is}
has index $12$ in $\PSL_{2}(\cO)$, as does $\gL$. 
\epf

The point  is that,
since
 $
\gL$ is now realized as an 
congruence
group,  its  elements can be parametrized,
giving an injection of affine space into the otherwise intractable thin Soddy group $\G$.
(In the Apollonian circle packing setting, the analogous idea 
was exploited extensively in, e.g., \cite{SarnakToLagarias, BourgainFuchs2011, BourgainKontorovich2014a}.)

\begin{prop}\label{prop:xiGgd}
For any $\g,
\gd
\in
\cO$ with 
\be\label{eq:gcd}
\g\equiv0(\mod\vr),\qquad
\lp\g,\gd\rp=\cO,
\ee 
there is an element
$$
\xi_{\g,\gd}:=
J^{-1}\cdot\rho\mattwo**{\g}\gd\cdot J
\in \Xi<\G_{1}<\G
,
$$
where 
$$
\xi_{\g,\gd}=
\bp
1&0&0&0&0\\
*&*&*&*&*\\
*&*&*&*&*\\
*&*&*&*&*\\
V&W&X&Y&Z
   \ep
,
$$
and
\beann
V&=& 
\frac{2}{3} \Re\left(\vr\gamma   \bar{\delta }\right)+|\gamma |^2+|\delta |^2-1
,\\
W&=&
   -
\frac{2}{3} \Re\left(
\gw\vr
\gamma  \bar{\delta }\right)
      ,\\
X&=&
-   \frac{2}{3} \Re\left(
   \bar\vr
   \gamma 
     \bar{\delta }\right)
,\\
Y&=&
   \frac{2}{3}
   \Re\left(  \vr \gamma  \bar{\delta }\right)
   +|\gamma |^2
,\\
Z&=&
\frac{2}{3} \Re\left(\vr\gamma    \bar{\delta }\right)
   +|\delta |^2
   .
\eeann
\end{prop}
\pf
This follows directly from \eqref{eq:Gam0Is}, \eqref{eq:Jcong}, \eqref{eq:Xidef}, \eqref{eq:XiIs} and  \eqref{eq:G1is}.
\epf

Recall  that  $\sO=\G\cdot\bb 
$ in \eqref{eq:cOAp} is the 
orbit
under the Soddy group $\G$
 of 
a
quintuple  $\bb
=(b_{1},\dots,b_{5})$ of bends.
According to \lemref{lem:loc}, there is an $\vep=\vep(\sP)\in\{\pm1\}$ so that every bend in $\sB$ is $\equiv0$ or $\vep(\mod 3)$.

Recalling that the set $\sB$ of bends contains sets of the form \eqref{eq:KwGv}, and setting $\bw=\bbe_{5}$, 
\propref{prop:xiGgd} immediately implies the following key
\begin{cor}\label{cor:fF}
Let $\bb\in\sO$ be a quintuple of bends,  and assume that
  $\g,\gd\in\cO$ satisfy \eqref{eq:gcd}. Then
the integer
\be\label{eq:fFbbDef}
\fF_{\bb}(\g,\gd)
:=
\<\bbe_{5},\xi_{\g,\gd}\cdot\bb\>
\ee
is in
the set $\sB$ of bends.
Setting
\be\label{eq:ffToF}
\ff_{\bb}
 \ := \
\fF_{\bb}+b_{1},
\ee
we have that $\ff_{\bb}$ is a homogeneous quaternary quadratic form 
given by:
\bea
\nonumber
\ff_{\bb}(\vr\g,\gd)
&=&
3A(\gamma _1^2+\gamma _1 \gamma _2  +\gamma _2^2)
+3
B(\gamma _1 \delta _1 
+\gamma _2 \delta _2)
   -3C\gamma _2 \delta _1
   \\
\label{eq:fvIs}
   &&
+3(B+C)\gamma _1 \delta
   _2 
   + D
   (\delta _1^2
   +\delta _1 \delta _2    +\delta _2^2)
.
\eea
Here the coefficients
$A,B,C,D$
are as in \eqref{eq:ABCD}, and $\g=\g_{1}+\g_{2}\gw$, $\gd=\gd_{1}+\gd_{2}\gw$ with $\g_{j},\gd_{j}\in\Z$.

Abusing notation, we write 
\be\label{eq:ffbx}
\ff_\bb(\bx)=\ff_\bb(\vr\g,\gd),
\ee 
where $\bx:=(\g_{1},\g_{2},\gd_{1},\gd_{2})$.
The (classically integral) symmetric matrix (that is, Hessian) corresponding to $\ff_{\bb}(\bx)$ is 
\be\label{eq:bbADef}
\bbA:=
\left(
\begin{array}{cccc}
6A& 3A &3B &3(B+C) \\
 3A & 6A &-3C&3B  \\
3B& -3C& 2 D &D \\
3(B+C)&3B&D & 2 D
\end{array}
\right)
,
\ee
so that $\ff_{\bb}(\vr\g,\gd)=\foh\bx\bbA\bx^{t}$. By \eqref{eq:cone}, the discriminant of $\ff_{\bb}$
is
\be\label{eq:disc}
\discr(\ff_{\bb})=|\bbA|=
9
\left(\foh Q(\bb)-3b_{1}^{2}
\right)^{2}
=
(3
b_{1})^{4}
.
\ee
Assume further that
\be\label{eq:order}
 b_{1}\le b_{2}\le b_{3}\le b_{4}\le b_{5}
,\quad\text{and}\quad
 b_{2}\ge0
.
\ee
Then the form
 $\ff_\bb$ is positive definite iff $b_1\neq0$ (otherwise it is positive semidefinite).
\end{cor} 
\pf
All direct computation. This should also elucidate the choice of the change of variables in \eqref{eq:ABCD}.
\epf
\begin{Def}
We say that $\fF_\bb$ ``$\cO$-primitively'' represents an integer $n$ if there exist $\g,\gd\in\cO$ satisfying \eqref{eq:gcd} so that $\fF_\bb(\g,\gd)=n$.
\end{Def}
We have thus shown that $\sB$ contains all the 
$\cO$-primitive
values of
the shifted quaternary quadratic form $\fF_\bb$. 
In the next section, we show that enough numbers are represented by such forms to produce a local-global principle in $\sB$.

\newpage

\section{Proof of  The Local-Global Theorem}

Recall from \lemref{lem:loc} that, to a primitive integral Soddy packing $\sP$, one assigns the number $\vep=\vep(\sP)\in\{\pm1\}$, so that 
every bend in $\sB=\sB(\sP)$ is congruent either $0$ or $\vep$ modulo $3$.
The analysis turns out to 
require that the odd primes dividing $b_1$ be  $\equiv1(\mod 3)$, so we first claim that this can always be arranged.

\begin{thm}\label{thm:makePrime}
If $\vep(\sP)=+1$, then there exists a (rational) prime 
\be\label{eq:fpDef}
\fp\equiv1(\mod 3)
\ee 
which
is a bend in $\sP$. If $\vep(\sP)=-1$, then $2\fp$ is a bend.
\end{thm}

Before giving the proof, we explain how this fact will be used.
By \corref{cor:fF}, we turn our attention to numbers $\cO$-primitively represented by $\fF_\bb$, as these are guaranteed to be in the bend set $\sB$.
It turns out that these are all $\equiv b_5(\mod 3)$, which is fine for our purposes, since we can make $b_5\equiv0$ or $\vep (\mod 3)$ by a choice of the quintuple $\bb=(b_j)$. 
Changing to the homogeneous form $\ff_\bb$ as in \eqref{eq:ffToF}, it will then suffice to show the following 
\begin{thm}\label{thm:FullLocGlob}
Assume that the quintuple $\bb$ has $b_1=\fp\equiv1(\mod 3)$ or $b_1=2\fp$, 
 and is ordered, that is, satisfies \eqref{eq:order}.
 Then every sufficiently large $n\equiv b_1+b_5(\mod 3)$ is $\cO$-primitively represented by $\ff_\bb$.
\end{thm}

Let
\be\label{eq:sRbbDef}
\sR_{\bb}(n):=\sum_{\g,\gd\in\cO\atop\lp\vr\g,\gd\rp=\cO
}\bo_{\{n=\ff_{\bb}(\vr\g,\gd)\}}
\ee
be the number of $\cO$-primitive representations of $n$ by $\ff_{\bb}$.
The study of this function will prove both theorems, with most of the tools going into the proof of the first also being useful for the second.
The 
key proposition which follows
is essentially
Kloosterman's method for representations by quaternary forms (as championed in this generality by Malyshev).

For an integer $m\ge1$ and a prime power $p^{a}$,
define the   $p$-adic local density $\gs_{p}(m;\bb)$ by
\be\label{eq:gsIs}
\gs_{p}(m;\bb)
:=
\lim_{a\to\infty}
\frac1{p^{3a}}
\#\{
\bx\in(\Z/p^{a}\Z)^{4}:
\ff_{\bb}(\bx)
\equiv m\ (\mod p^{a})
\}
,
\ee
where we have used the convention \eqref{eq:ffbx}.

\begin{prop}\label{prop:Kloo}
 If $n\equiv b_1+b_5(\mod 3)$ and $(b_1,3)=1$, then
\be\label{eq:propKloo}
\sR_\bb(n) \ =\ 
n\,
{
\pi
^{2}
\over 
9 b_{1}^{2}
}
\fS_0(n;\bb)\fS_1(n;\bb)\fS_2(n;\bb)
+ 
O_{\bb,\gep}\left(
n^{3/4+\gep}
\right)
,
\ee
with an effective implied constant.
Here
\be\label{eq:fSnbbDef}
\fS_0(n;\bb) \ :=\
\prod_{p}
\gs_p(n;\bb),
\ee
$$
\quad
\fS_1(n;\bb) \ :=\
\prod_{p\equiv1(3)\atop p\mid n}
\gs^{(1)}_p(n;\bb),
\qquad
\fS_2(n;\bb) \ :=\
\prod_{p\equiv2(3)\atop p^{2}\mid n}
\gs^{(2)}_p(n;\bb),
$$
where
\be\label{eq:Type1Def}
\gs_p^{(1)}(n;\bb)\ :=\
\left(
1
-
{2\over p }{
\gs_{p}\left({n\over p};\bb\right)
\over
\gs_{p}\left({n};\bb\right)
}
+\bo_{\{p^{2}\mid n\}}
{1\over p^{2}}{
\gs_{p}\left({n\over p^{2}};\bb\right)
\over
\gs_{p}\left({n};\bb\right)
}
\right)
,
\ee
and
\be\label{eq:Type2Def}
\gs_p^{(2)}(n;\bb)\ :=\
\left(
1-{1\over p^{2}}{
\gs_{p}\left({n\over p^{2}};\bb\right)
\over
\gs_{p}\left({n};\bb\right)
}
\right)
.
\ee
\end{prop}

We will call the terms arising in $\fS_j$ ``Type $j$'', and refer to primes $p$ as ``Good'' or ``Bad'' depending on whether $(p, 2\cdot3\cdot b_1)=1$ or not.
While the proof largely uses standard techniques, a few of the  manipulations are somewhat delicate, so
we give the details. 

\pf
Recall that the Dedekind zeta function of $K$ is
$$
\gz_{K}(s):=\sum_{\fm}{1\over \N\fm^{s}}=\prod_{\fp}\left(1-{1\over\N\fp^{s}}\right)^{-1},
$$
where $\N$ is the norm, the sum is over non-zero integral ideals $\fm$ of $K$, and the product is over prime ideals $\fp$. We define the $K$-M\"obius function $\mu_{K}$ via
$$
{1\over \gz_{K}(s)} =\prod_{\fp}\left(1-{1\over\N\fp^{s}}\right)=\sum_{\fm}{\mu_{K}(\fm)\over \N\fm^{s}}.
$$
Thus $\mu_{K}$ is multiplicative, supported on non-zero, square-free integral ideals, and takes the value $-1$ on prime ideals. M\"obius inversion now reads:
$$
\sum_{\fd\supset\fm}\mu_{K}(\fd)=\twocase{}{1}{if $\fm=\cO$,}{0}{otherwise.}
$$

M\"obius inversion works on the level of ideals, but $\ff_{\bb}$ in \eqref{eq:sRbbDef} is  a function on elements of $\cO$, not ideals (i.e. it is {\it not} invariant under units in each variable $\g,$ $\gd$ separately). So we will have to pass from ideals to elements, and back again. 
Begin by writing
\beann
\sR_{\bb}(n)
&=&\sum_{\g,\gd\in\cO}\bo_{\{n=\ff_{\bb}(\vr\g,\gd)\}}\sum_{\fd\supset\lp\vr\g,\gd\rp}\mu_{K}(\fd)
\\
&=&
\sum_{\fd}\mu_{K}(\fd)
\sum_{\g,\gd\in\cO\atop\lp\vr\g\rp\subset\fd,\lp\gd\rp\subset\fd}\bo_{\{n=\ff_{\bb}(\vr\g,\gd)\}}
.
\eeann

The field $K$ is a principal ideal domain with a finite group of units, $|\cO^{\times}|=6$, so the non-zero integral ideals of $K$ are in 1-to-6 correspondence with non-zero elements of $\cO$. So we can write
$\fd=\lp\eta\rp$ with $\eta\in\cO\setminus0$, whence
\beann
\sR_{\bb}(n)
&=&
\frac1{|\cO^{\times}|}
\sum_{\eta\in\cO}\mu_{K}(\lp\eta\rp)
\sum_{\g,\gd\in\cO\atop\vr\g\equiv0(\mod\eta),\gd\equiv0(\mod\eta)}\bo_{\{n=\ff_{\bb}(\vr\g,\gd)\}}
.
\eeann

Now comes a little trick which will allow us to replace $\vr\g\equiv0(\eta)$ by just $\g\equiv0(\eta)$. 
Indeed, an easy calculation shows that
\be\label{eq:ffHom}
\ff_{\bb}(\eta\g',\eta\gd')=
\N\eta
\cdot
\ff_{\bb}(\g',\gd').
\ee
So $n=\ff_\bb(\vr\g,\gd)$, together with $\vr\g,\gd\equiv0(\mod \eta)$, implies that $\N\eta$ divides $n$. 
But  $b_1\equiv\vep(\mod 3)$, $b_5\equiv0$ or $\vep(\mod 3)$, and $n\equiv b_1+b_5(\mod 3)$ together imply that $n\equiv\vep$ or $2\vep(\mod 3)$. In particular, $(n,3)=1$, hence $(\N\eta,3)=1$, so $\vr$ is coprime to $\eta$. Now we have:
\beann
\sR_{\bb}(n)
&=&
\frac1{|\cO^{\times}|}
\sum_{\eta\in\cO\atop
\N\eta\mid n
}\mu_{K}(\lp\eta\rp)
\sum_{\g,\gd\in\cO
}\bo_{\{{n\over\N\eta}=\ff_{\bb}(\vr\g,\gd)\}}
.
\eeann

Having freed the variables $\g,\gd$, we may return to ideals, and use the convention \eqref{eq:ffbx} to write
\be\label{eq:sRtocR}
\sR_{\bb}(n)
=
\sum_{
\N\fd\mid n}\mu_{K}(\fd)
\cR_{\bb}\left({n\over\N\fd}\right)
,
\ee
where
$$
\cR_{\bb}(m):=
\sum_{\bx
\in\Z^{4}}
\bo_{\{m=\ff_{\bb}(\bx)
\}}
$$
is now a classical representation quantity. 

Combining \eqref{eq:disc} with \cite[(11.57), (11.62), (11.19)]{Iwaniec1997book},
we have, for any $\gep>0$ (not to be confused with $\vep(\sP)\in\{\pm1\}$),
\be\label{eq:Kloo}
\cR_{\bb}(m)
=
{
\pi
^{2}
\over 
9 b_{1}^{2}
}\,
m\,
\fS_0(m,\bb)
\ +\ 
O_{\bb,\gep}(m^{3/4+\gep})
,
\ee
where the singular series
$
\fS_0(m,\bb)
$
is as in \eqref{eq:fSnbbDef}.
The implied constant is effective.
 Inserting \eqref{eq:Kloo} into \eqref{eq:sRtocR} gives
\bea
\label{eq:sRbbKloo}
\sR_{\bb}(n)
&=&
n\,
{
\pi
^{2}
\over 
9 b_{1}^{2}
}
\sum_{\N\fd\mid n}{\mu_{K}(\fd)\over\N\fd}
\fS_0\left({n\over\N\fd},\bb\right)
\\
\nonumber
&& \hskip1in+ 
O_{\bb,\gep}\left(
n^{3/4+\gep}
\sum_{\N\fd\mid n}1
\right)
.
\eea
We clearly have $\sum_{\N\fd\mid n}1\ll_{\gep} n^{\gep}$, so the error term is as claimed in \eqref{eq:propKloo}. It remains to control the local densities. 

Recall that $\fd$ is a square-free ideal.
Let $p$ be a rational prime dividing $\N\fd$; then $p\neq3$. If $p\equiv2(3)$ is inert, then $\N (p)=p^{2}$ and we can write $\fd=\lp p\rp\fd'$, where $(\N\fd',p)=1$; thus $\ord_{p}(\N\fd)=2$. 
If $p\equiv1(3)$ splits in $\cO$ as $\lp p\rp=\pi\bar\pi$, then we have $\ord_{p}(\N\fd)=2$ or $1$, depending on whether both $\pi$ and $\bar\pi$ divide $\fd$ or just one.  
Either way, we can write 
$$
\fd=\fp_{0}\fd'\quad\text{ with }\quad(\N\fd',p)=1.
$$ 
Extend this notation to rational primes $p$ not dividing $\N\fd$ by setting $\fp_{0}=\cO$ and $\fd'=\fd$. We claim that:
\be\label{eq:densitiesEq}
\gs_{p}\left({n\over \N\fd};\bb\right)
=
\gs_{p}\left({n\over \N\fp_{0}};\bb\right)
.
\ee
Indeed,
let $\fd'=(\eta')$.
Applying \eqref{eq:ffHom} in reverse together with \eqref{eq:gsIs}, 
we see that
the density 
$
\gs_{p}\left({n\over \N\fd};\bb\right)
$
is counting the number of solutions  to 
$$
\ff_{\bb}(\eta' \vr(\g_{1}+\g_{2}\gw),\eta'(\gd_{1}+\gd_{2}\gw))\equiv {n\over\N\fp_{0}}\qquad(\mod p^a).
$$
The linear map 
$$
(\vr(\g_{1}+\g_{2}\gw),\gd_{1}+\gd_2\gw)\mapsto (\eta'\vr (\g_{1}+\g_{2}\gw),\eta'(\gd_{1}+\gd_{2}\gw))
$$ 
has determinant $(\N\fd')^{2}$, and hence is invertible since $(p,\N\fd')=1$.
Thus the two densities agree and we have proved \eqref{eq:densitiesEq}.

In particular, we have
\be\label{eq:hDef}
\sum_{\N\fd\mid n}{\mu_{K}(\fd)\over\N\fd}
\fS_0\left({n\over\N\fd},\bb\right)
=
\fS_0\left({n},\bb\right)
\sum_{\N\fd\mid n}{\mu_{K}(\fd)\over\N\fd}
h(\fd)
,
\ee
where
$$
h(\fd):=
\prod_{p\mid \N\fd\atop p\equiv2(3)}
{
\gs_{p}\left({n\over p^{2}};\bb\right)
\over
\gs_{p}\left({n};\bb\right)
}
\prod_{p\mid \N\fd\atop p\equiv1(3)}
{
\gs_{p}\left({n\over p^{\ord_{p}(\N\fd)}};\bb\right)
\over
\gs_{p}\left({n};\bb\right)
}
,
$$
assuming the denominators do not vanish. Note that this function is {\it not} multiplicative
on ideals.
Nevertheless, 
we
do 
have a factorization with respect to rational primes of the following form:
\bea
\nonumber
&&
\hskip-.5in
\sum_{\N\fd\mid n}{\mu_{K}(\fd)\over\N\fd}
h(\fd)
=
\prod_{p\equiv2(3)\atop p^{2}\mid n}
\left(
1+{\mu_{K}(\lp p\rp)\over \N (p)}{
\gs_{p}\left({n\over p^{2}};\bb\right)
\over
\gs_{p}\left({n};\bb\right)
}
\right)
\\
\nonumber
&&
\hskip-.5in
\times
\prod_{p\equiv1(3)\atop p\mid n,(p)=\pi\bar\pi}
\left(
1
+
{\mu_{K}(\pi)\over \N \pi}{
\gs_{p}\left({n\over p};\bb\right)
\over
\gs_{p}\left({n};\bb\right)
}
+
{\mu_{K}(\bar\pi)\over \N \bar\pi}{
\gs_{p}\left({n\over p};\bb\right)
\over
\gs_{p}\left({n};\bb\right)
}
+\bo_{\{p^{2}\mid n\}}
{\mu_{K}(\pi\bar\pi)\over \N(\pi \bar\pi)}{
\gs_{p}\left({n\over p^{2}};\bb\right)
\over
\gs_{p}\left({n};\bb\right)
}
\right)
\\
\nonumber
&&
=
\prod_{p\equiv2(3)\atop p^{2}\mid n}
\left(
1-{1\over p^{2}}{
\gs_{p}\left({n\over p^{2}};\bb\right)
\over
\gs_{p}\left({n};\bb\right)
}
\right)
\prod_{p\equiv1(3)\atop p\mid n}
\left(
1
-
{2\over p }{
\gs_{p}\left({n\over p};\bb\right)
\over
\gs_{p}\left({n};\bb\right)
}
+\bo_{\{p^{2}\mid n\}}
{1\over p^{2}}{
\gs_{p}\left({n\over p^{2}};\bb\right)
\over
\gs_{p}\left({n};\bb\right)
}
\right)
.
\\
\label{eq:hIs}
\eea
The problem has now returned back to the rational integers $\Z$ after the detour through the Eisenstein ones $\cO$ . 
Inserting \eqref{eq:hIs} and \eqref{eq:hDef} into \eqref{eq:sRbbKloo} gives \eqref{eq:propKloo}, as claimed.
\epf

It remains to analyze the Type 0, 1, and  2 factors for both Good and Bad primes. First the Good.

\begin{lem}\label{lem:Good}
Let $n$ and $b_1$ as in \propref{prop:Kloo}, and 
assume that $p$ is Good, that is, $(p,2\cdot3\cdot b_1)=1$. Then 
\be\label{eq:Type0}
\gs_p(n;\bb)=\twocase{}
{1+\frac1p+O({1\over p^2})}
{if $p\mid n$,}
{1-{1\over p^2}}
{otherwise,}
\ee
$$
\gs^{(1)}_p(n;\bb)=\twocase{}
{1+O(\frac1p)}
{if $p\mid n$,}
{1}
{otherwise,}
$$
$$
\gs^{(2)}_p(n;\bb)=
{1+O\left({1\over p^2}\right)}
,
$$
and none of these local factors vanish.
\end{lem}
\pf
Write $p^{k}\| n$.
We first handle $\gs_p(n;\bb)$.
Apply \cite[(11.69), (11.70), (11.72)]{Iwaniec1997book}, giving
$$
\gs_{p}\left({n};\bb\right)
=
{\left(1-{\chi(p)\over p^{2}}\right)\left(1-{\chi(p^{k+1})\over p^{k+1}}\right)
\over
\left(1-{\chi(p)\over p}\right)
}
.
$$
Here $\chi(m):=\left({|\bbA|\over m} \right)$ is the quadratic character modulo the discriminant $|\bbA|=(3 b_{1})^{4}$ of $\ff_\bb$, cf. \eqref{eq:disc}. Since the latter is a
 square and  $(p,|\bbA|)=1$, we have that $\chi(p)=\chi(p^k)=1$;  
hence
\be\label{eq:gspIs}
\gs_{p}\left({n};\bb\right)
=
{\left(1+{1\over p}\right)\left(1-{1\over p^{k+1}}\right)
}
.
\ee
This clearly never vanishes, and \eqref{eq:Type0} is readily verified.

Next we deal with $\gs_p^{(2)}$. Then $p\equiv2(3)$ and $p^{2}\mid n$, that is, $p^{k}\|n$ with $k\ge2$. Inserting \eqref{eq:gspIs} into \eqref{eq:Type2Def} gives
$$
\gs_p^{(2)}(n;\bb)=
1-{1\over p^{2}}{
\gs_{p}\left({n\over p^{2}};\bb\right)
\over
\gs_{p}\left({n};\bb\right)
}
=
1-{1\over p^{2}}{
\left(1-{1\over p^{k-1}}\right)
\over
\left(1-{1\over p^{k+1}}\right)
}
=
1-
{
p^{k-1}-1
\over
p^{k+1}-1
}
.$$
This factor clearly never vanishes, and for $p$ large is of size
$
1
+O\left({1\over p^{2}}\right)
,
$
so
is harmless. 

Finally we handle $\gs_p^{(1)}$.
Here $p\equiv1(3)$ and there are two cases depending on whether  $k=1$, or $k\ge1$. If $k=1$, then the factor is
$$
\gs_p^{(1)}(n;\bb)=
1
-
{2\over p }{
\gs_{p}\left({n\over p};\bb\right)
\over
\gs_{p}\left({n};\bb\right)
}
=
1
-
{2\over p }{
{\left(1-{1\over p}\right)
}
\over
{\left(1-{1\over p^{2}}\right)
}
}
=
1
-
{
2
\over
p+1
}
,
$$
which doesn't vanish.
If $k\ge2$, then the factor is
\beann
&&
1
-
{2\over p }{
\gs_{p}\left({n\over p};\bb\right)
\over
\gs_{p}\left({n};\bb\right)
}
+
{1\over p^{2}}{
\gs_{p}\left({n\over p^{2}};\bb\right)
\over
\gs_{p}\left({n};\bb\right)
}
=
1
-
{2\over p }{
\left(1-{1\over p^{k}}\right)
\over
\left(1-{1\over p^{k+1}}\right)
}
+
{1\over p^{2}}{
\left(1-{1\over p^{k-1}}\right)
\over
\left(1-{1\over p^{k+1}}\right)
}
\\
&&
\hskip.5in
=
1
-
{
 2
p^{k}
-p^{k-1}
-1
\over
p^{k+1}-1
}
.
\eeann
This again does not vanish, and is asymptotically of size $1+O(\frac1p)$. This completes the proof.
\epf

To deal with the Bad primes, we first record Hensel's Lemma. Recall from \eqref{eq:bbADef} that $\bbA$ is the Hessian of $\ff_\bb$.

\begin{lem}[Hensel's Lemma]\label{lem:Hensel}
Assume that
$$
\ff_\bb(\bx)
\equiv n(\mod p^k)
$$
for $\bx\in(\Z/p^k\Z)^4$ with $\bx\bbA\neq0(\mod p)$. Then
the set of ``lifts'' $\by\in(\Z/p^{k+1}\Z)^4$ with $\by\equiv\bx(\mod p^k)$ having
\be\label{eq:ffby}
\ff_\bb(\by)\equiv n(\mod p^{k+1})
,
\ee
has cardinality exactly $p^3$.

If on the other hand $\bx\bbA\equiv0(\mod p)$, then the number of lifts is either $p^4$ or $0$, depending on whether $\ff_\bb(\bx)\equiv n(\mod p^{k+1})$ or not.
\end{lem}
\pf
Write $\by=\bx+p^k\bba$ with $\bba\in(\Z/p\Z)^4$. 
The equation
$$
\ff_\bb(\bx+p^k\bba)\ = \ \ff_\bb(\bx)+p^k\bx\bbA\bba^t+p^{2k}\ff_\bb(\bba)
$$
is valid in the integers, and hence also valid mod $p^{k+1}$, even for $p=2$. Write 
$$
n-\ff_\bb(\bx)\equiv p^k m\qquad(\mod p^{k+1})
.
$$ 
Then the equation \eqref{eq:ffby} becomes
\be\label{eq:Hensel}
m=
\bx\bbA\bba^t
(\mod p
)
.
\ee
If $\bx\bbA$ is not the zero vector mod $p$, then there are exactly $p^3$ solutions for $\bba$, and hence for $\by$, as claimed.

If, on the other hand, $\bx\bbA\equiv0(\mod p)$, then \eqref{eq:Hensel} has either $p^4$ or $0$ solutions, depending on whether $m\equiv0$ or not, that is, whether $\ff_\bb(\bx)\equiv n(\mod p^{k+1})$ or not.
\epf

This is already sufficient to deal conclusively with the crucial prime $p=3$, for which there is only Type 0, and the only relevant $n$'s are those coprime to $3$. 
\begin{lem}\label{lem:p3}
Let $n$ and $b_1$ as in \propref{prop:Kloo}.
 Then
$$
\gs_3(n;\bb) 
\ \gg \ 1.
$$
\end{lem}
\pf
Reducing \eqref{eq:fvIs} mod $3$ shows that $\ff_{\bb}(\vr\g,\gd)\equiv( b_{1}+ b_{5})\N(\gd)$, where $\N(\gd)=\gd_{1}^{2}+\gd_{1}\gd_{2}+\gd_{2}^{2}$. As $\vr$ and $\gd$ are coprime, $\N(\gd)\equiv1(\mod 3)$. Hence
$\ff_{\bb}(\vr\g,\gd)$ is always
$$
\equiv b_{1}+ b_{5}(\mod3)
.
$$
Having assumed that $b_1\equiv\vep(\mod 3)$, we have that $D=b_1+b_5\equiv\vep$ or $2\vep(\mod 3)$, in either case this is coprime to $3$. Hence Hensel's lemma applies and solutions can be lifted to the $3$-adic integers $\Z_3$.
\epf

Next we record that for Bad primes $p\neq3$, the Hessian $\bbA$ cannot vanish completely.
\begin{lem}\label{lem:HessModP}
Assume $p\mid 2b_1$ and $p\neq3$. Then $\bbA$ is not identically zero mod $p$.
\end{lem}
\pf
For $p=2$, this is a direct calculation. Indeed, if $\bbA\equiv0(\mod p)$, then 
\be\label{eq:ABCBCD}
A\equiv B \equiv C\equiv D\equiv B+C,
\ee
which forces the $b_j$'s to either be all $0$ or all $1$. The former is impossible by the primitivity of $\bb$.
The latter is also impossible from looking at the cone \eqref{eq:cone} mod $4$.

If instead $p\mid b_1$, then \eqref{eq:ABCBCD} forces the $b_j$ to all be $\equiv0(\mod p)$, which again is impossible by primitivity.
\epf

By \lemref{lem:HessModP}, the Hessian $\bbA$ has a non-zero entry; assume that $A\neq0(\mod p)$, the other cases being similar.
For Bad primes $p\ge5$, that is, those diving $b_1$, the following is a convenient normal basis for studying the quadratic space of $\ff_\bb$.

\begin{lem}\label{lem:A0p}
Assume $p\mid b_1$, $p\ge5$, and $A\neq0(\mod p)$. Then
 the following vectors 
$$
\bu_1=(1,0,0,0),\
\bu_2=(1,-2,0,0),
$$
$$
\bu_3=(-B-2C,-B+C,0,3A),\
\bu_4=(-2B-C,B+2C,3A,0),
$$
form a basis for $\F_p^4$ which
 is normal, 
that is, $\bu_j\bbA\bu_k\equiv0$ if $j\neq k$.

Moreover, 
\be\label{eq:fbujs}
\ff_\bb(\bu_1)=3A,\ 
\ff_\bb(\bu_2)=9A,\ 
\ff_\bb(\bu_3)=\ff_\bb(\bu_4)=-9A\cdot F(\bba), 
\ee
where $F(\bba)$ is given in \eqref{eq:Fis}.
Hence writing any $\bx\in\F_p^4$ as 
\be\label{eq:bxabcd}
\bx=a\bu_1+b\bu_2+c\bu_3+d\bu_4,
\ee 
we  have that
\be\label{eq:ffbxaa3bb}
\ff_\bb(\bx)\ \equiv \ 3Aa^2+9Ab^2\ \equiv \ 3A(a^2+3b^2)\quad(\mod p).
\ee
\end{lem}
\pf
By 
\eqref{eq:Ftogk1}, 
$$
F(\bba)\ = \ B^{2}+BC+C^{2}-AD\ \equiv \ 0(\mod p)
,
$$
whence 
$\bu_3$ and $\bu_4$ 
are null vectors for $\bbA$, that is, 
$$
\bu_3\bbA \equiv\bu_4\bbA \equiv0(\mod p).
$$
The rest is readily verified by computation.
\epf

The appearance of the binary form $a^2+3b^2$ in \eqref{eq:ffbxaa3bb} 
explains
why we want $b_1$ to contain only primes $p\equiv1(\mod 3)$ in \thmref{thm:makePrime};
indeed, if there are Bad primes $p\equiv2(\mod 3)$, then there {\it can} be further local obstructions mod $p^2$, and $\gs_p$ can vanish!
But first, we are now in position to give a

\subsection{Proof of \thmref{thm:makePrime}}\

Assume first that $\vep=\vep(\sP)=+1$.
By \lemref{lem:bend24mod6}, we may arrange $\bb$ so that
$$
b_1\equiv3+\vep=4
\qquad(\mod 6).
$$
 In particular, $b_1$ is even, and $b_1\equiv1(\mod 3)$. 
We may also assume that $b_5\equiv\vep\equiv1(\mod 3)$, as can be arranged by \lemref{lem:loc}.
We first claim that $\ff_\bb$, the homogeneous form, $\cO$-primitively represents every sufficiently large 
\be\label{eq:n2vep}
n\equiv 1+4b_1\qquad(\mod 6b_1).
\ee
Indeed, in this progression, 
$$
n\equiv 2\equiv2\vep\equiv b_1+b_5\ (\mod 3),
$$ 
so the conditions of \propref{prop:Kloo} is satisfied. Moreover, $n$
is coprime to $2 b_1$, 
so there are no Bad factors of Type 1 or 2, and Hensel's lemma, together with \lemref{lem:HessModP}, allows us to control the local densities at $2$ and at primes dividing $b_1$. Then \eqref{eq:propKloo} is a true asymptotic, giving the claim.

Returning to the shifted quaternary form $\fF_\bb$ in \eqref{eq:fFbbDef}, we have from \eqref{eq:ffToF} that
every sufficiently large value of 
\be\label{eq:nMinB1}
n-b_1\equiv 1+3b_1(\mod 6b_1)
\ee
is $\cO$-primitively represented by $\fF_\bb$, and hence appears in the set $\sB$ of bends. 
Such numbers are all $\equiv1(\mod 6)$, and
this arithmetic progression has coprime modulus and shift (since $b_1$ is even), whence Dirichlet's theorem applies, showing that $\sB$ contains a prime $\fp\equiv1(\mod 6)$. This
of course is equivalent to $\fp\equiv1(\mod 3)$. 

The argument for the case $\vep(\sP)=-1$ is similar, so we omit it. This completes the proof of \thmref{thm:makePrime}. 

\subsection{Proof of \thmref{thm:FullLocGlob}}\

Now we assume that $b_1=\fp$ or $b_1=2\fp$
as in \thmref{thm:makePrime}; clearly then $b_1\equiv\vep(\mod 3)$.
As before, if $n\equiv b_1+b_5(\mod 3)$ is coprime to $2b_1$ and sufficiently large, then it is $\cO$-primitively represented by $\ff_\bb$. So if $n$ is even, we need to handle the $2$-adic densities, both of Type 0 and Type 2, and when $n\equiv0(\mod \fp)$, we need control on  the $\fp$-adic local factors of Type 0 and Type 1.

We begin by recording the following
\begin{lem}\label{lem:429}
If $n\equiv0(\mod \fp),$ then
$$
\#\{\bx\in\F_\fp^4\ : \ \ff_\bb(\bx)\equiv n ,\ \bx\bbA\neq0\ (\mod \fp)\}
 \ =\
 2(\fp-1)\fp^2,
$$
and
$$
\#\{\bx\in\F_\fp^4\ : \ \ff_\bb(\bx)\equiv n ,\ \bx\bbA\equiv0\ (\mod \fp)\}
 \ =\
\fp^2.
$$
If $(n,\fp)=1$, then
$$
\#\{\bx\in\F_\fp^4\ : \ \ff_\bb(\bx)\equiv n\ (\mod \fp)\}
 \ =\
 (\fp-1)\fp^2.
$$
\end{lem}
\pf
This follows easily from \lemref{lem:A0p}. Indeed, assume that $A\neq0(\mod \fp)$, and first check the case $n\equiv0(\mod \fp)$. Then by \eqref{eq:ffbxaa3bb}, we
need (since the values $c$ and $d$ in \eqref{eq:bxabcd} are completely free) to count the number of $a^2+3b^2\equiv0(\mod \fp)$. Since $\fp\equiv1(\mod 3)$, there are $2(\fp-1)$ such non-trivial solutions%
\footnote{In the language of \cite[Ch. 2.2]{Cassels1978}, the span of $\bu_1$ and $\bu_2$ in \lemref{lem:A0p} is a regular, isotropic subspace of $\F_\fp^4$ when $\fp\equiv1(\mod 3)$.}%
, plus one trivial, $(a,b)=(0,0)$. For any of these, $\bx\bbA\equiv0$ if and only if $(a,b)$ is trivial, so the total number of solutions is as claimed.

If $(n,\fp)=1$, then the number of solutions, say $\cN$, is independent of $n$. By the counts for $\ff_\bb(\bx)\equiv0$ above, we then have that
$$
(\fp-1)\cN+\fp^2+2(\fp-1)\fp^2 \ =\  \fp^4,
$$
since there are $\fp^4$ total choices for $\bx$. Solving for $\cN$ gives the claim.
\epf
Lifting these solutions by Hensel's lemma, we completely control the Type 0 factors, as follows.
\begin{lem}
If $(n,\fp)=1$, then
\be\label{eq:BadLocDensP0np}
\gs_\fp(n;\bb) \ = \ \left(1-{1\over \fp}\right).
\ee
If $\fp\| n$, then 
$$
\#\{\bx\in(\Z/\fp^2\Z)^4\ : \ \ff_\bb(\bx)\equiv n  (\mod \fp^2)\} \ = \ 2 (\fp - 1) \fp^5,
$$
whence
\be\label{eq:BadLocDensP0}
\gs_\fp(n;\bb) \ = \ 2\left(1-{1\over \fp}\right).
\ee
If $\fp^2\mid n$, then
$$
\#\{\bx\in(\Z/\fp^2\Z)^4\ : \ \ff_\bb(\bx)\equiv n  (\mod \fp^2)\} \ = \ 2 (\fp - 1) \fp^5+\fp^6,
$$
and
\be\label{eq:BadLocLower}
\gs_\fp(n;\bb) \ \ge \ 2\left(1-{1\over \fp}\right).
\ee
\end{lem}
\pf
If $(n,\fp)=1$, then the last statement of \lemref{lem:429} applies, and can be lifted by Hensel's Lemma, giving
\eqref{eq:BadLocDensP0np}.

Next consider the case $\fp\|n$. By \lemref{lem:429}, 
there are $2(\fp-1)\fp^2$ ``non-trivial'' solutions mod $\fp$ (i.e., those with $\bx\bbA\neq0$), and by Hensel's Lemma, these each lift to 
$\fp^3$ solutions mod $\fp^2$. We claim that the trivial mod $\fp$ solutions (those with $\bx\bbA\equiv0$) have no lifts mod $\fp^2$. Indeed, 
\eqref{eq:fbujs} and \eqref{eq:Ftogk1} imply that $\ff_\bb(\bu_3)\equiv\ff_\bb(\bu_4)\equiv0(\mod \fp^2)$, and  $n/\fp$ is coprime to $\fp$, so the trivial solutions do not lift.
 This gives the asserted count, and also \eqref{eq:BadLocDensP0} by iterating Hensel's Lemma.

If $\fp^2\mid n$, then the non-trivial mod $\fp$ solutions 
still each lift to $\fp^3$ solutions mod $\fp^2$. But now the trivial mod $\fp$ solutions 
also lift, and each has $\fp^4$ lifts, since $\bba$ in \eqref{eq:Hensel} is completely free. The lower bound \eqref{eq:BadLocLower} comes from lifting just the non-trivial solutions.
\epf

Thus the Type 0 local density is controlled. We can also now handle the Type 1 local density. 
\begin{lem}
If  $\fp\|n$, then
$$
\gs_\fp^{(1)}(n;\bb) \ = \ 1-\frac1\fp.
$$
If $p^2\mid n$, then
$$
\gs_\fp^{(1)}(n;\bb) \ \ge \ 1-\frac1\fp.
$$
\end{lem}
\pf
If $p\|n$, then the claim follows trivially on combining 
\eqref{eq:BadLocDensP0np} and \eqref{eq:BadLocDensP0} into \eqref{eq:Type1Def}, where there is no third term.

If $p^k\| n$ with $k\ge2$, then there is a third term in \eqref{eq:Type1Def}, but we can drop it by positivity (since we're only claiming a lower bound). In the expression \eqref{eq:gsIs} for the local density $\gs_\fp$, the limit stabilizes as soon as $a>k+\ord_\fp(|\bbA|)$, so we can take $a=k+3$, since $\ord_\fp(|\bbA|)=2$. Setting
$$
\cN_\fb(n;p^a):=
\#\{
\bx\in(\Z/p^{a}\Z)^{4}:
\ff_{\bb}(\bx)
\equiv m\ (\mod p^{a})
\}
,
$$
we see that, since $p^k\|n$,
$$
\cN_\fb(n;p^{k+3}) \ \ge \
\cN_\fb(n/p;p^{k+3}),
$$
since the former may have more ``trivial'' lifts. Hence $\gs_\fp(n;\fb)\ge \gs_\fp(n/p;\fb)$, from which the claim follows.
\epf

This completes our analysis for the special Bad prime $p=\fp$. It remains to handle $p=2$. 

\begin{lem}
For $p=2$,
$$
\gs_2(n;\bb)\ \gg \ 1, \qquad \gs_2^{(2)}(n;\bb)\ \gg\ 1.
$$
\end{lem}
\pf
Assume first that $b_1,$ $b_2$, and $b_3$ are odd, and that $b_4$ is even, the other cases being similar. 
Reducing
\eqref{eq:ABCD} mod $2$ gives
$$
A\equiv1,\ 
B\equiv b_5,\
C\equiv b_5,\
B+C\equiv0,\
D\equiv1+b_5.
$$
For either possible value of $b_5$, there are six $\bx\in\F_2^4$ with $\ff_\bb(\bx)\equiv1$ and the other ten have $\ff_\bb(\bx)\equiv0$. 
One of the ten is of course the zero vector, and the remaining nine all have $\bx\bbA\neq0(\mod 2)$. Hence they lift by Hensel's lemma, giving  control on both $\gs_2$ and $\gs_2^{(2)}$. 
\epf

\subsection{Proof of \thmref{thm:main}}\

Now we put evertything together. By \thmref{thm:makePrime}, we take $b_1=\fp$ or $2\fp$, the arrange for the ordering \eqref{eq:order} to be satisfied. By \thmref{thm:FullLocGlob}, $\ff_\bb$ then $\cO$-primitively represents every large $n\equiv b_1+b_5(\mod 3)$. Hence $\fF_\bb=\ff_\bb-b_1$, the shifted form, $\cO$-primitively represents every large $n\equiv b_5(\mod 3)$,
and by \corref{cor:fF}, these numbers are all in $\sB$. Since we can make $b_5\equiv0$ or $\vep$ mod $3$, this covers all the local obstructions in \lemref{lem:loc}. In particular, they are {\it a posteriori} all the local obstructions. This completes the proof of the Local-Global Theorem.

\subsection{Explicit Example}\

We illustrate here the procedure described above for the example of the packing $\sP_0$ having ``root''  quintuple $\bb_0=(-11,21,25,27,28)$ as in \eqref{eq:bv0Is}. In this case, $\vep(\sP)=+1$, but $b_1=-11$ has prime factors (namely, $11$) which are not $\equiv 1(\mod 3)$, so we cannot apply \thmref{thm:FullLocGlob} directly.
Following 
 the proof of \thmref{thm:makePrime}, 
we first arrange for $b_1$ to be $\equiv4(\mod 6)$ and $b_5\equiv 1(\mod 3)$, by reordering $\bb_0$ to $\bb_1=(28,21,25,27,-11)$. This does not satisfy \eqref{eq:order}, so we apply $M_4M_5M_4M_3M_2M_5$ in \eqref{eq:Mj} to $\bb_1$, giving
$$
\bb_2=(28, 171, 313, 912, 997).
$$
(Note that at no point are we changing the bends appearing in $\sP_0$, and each quintuple still represents the bends of five mutually tangent spheres. We also only apply even length words in $M_j$, $j=2,\dots, 5$, so are staying within $\Xi$ in \eqref{eq:XiIs}.) 
Now we have $b_1\equiv4(\mod 6)$ and $b_5\equiv1(\mod 3)$, so can argue as in \eqref{eq:nMinB1} to show that  the set $\sB$ of bends contains all sufficiently large values of the progression $85(\mod 168)$. The smallest of these, $\fp=421$, turns out to already be in $\sB$;
in fact, applying $M_5M_3M_4M_3M_5M_4$ to $\bb_2$, and reordering to make $b_1=\fp$  gives 
$$
\bb_3=({ 421, 25, 28, 171, 309 })
.
$$ 
Now apply $\G$ some more to correct the ordering, 
$$
\bb_4=
M_5M_4M_3M_2\cdot\bb_3=({421, 904, 1777, 3240, 6033}).
$$ 
We are finally in position to apply \thmref{thm:FullLocGlob}; then every sufficiently large number 
$$
n\equiv b_5\equiv0(\mod 3)
$$ 
is $\cO$-primitively represented by the shifted form $\fF_\bb$, and hence appears in $\sB$ by \corref{cor:fF}. 
Next we apply $M_5.M_3$ to $\bb_4$ and reorder to obtain
$
\bb_5 = ({421, 904, 3240, 7353, 8821}).
$
This has $b_5\equiv1(\mod 3)$, and hence all large numbers $\equiv1(\mod 3)$ also appear in $\sB$, as claimed.

\newpage
\bibliographystyle{alpha}

\bibliography{AKbibliog}

\end{document}